 \documentclass[final,3p,times]{elsarticle}

\usepackage{color}
\usepackage{multirow}
\usepackage{longtable}
\usepackage{footnote}

\usepackage{amssymb}

\journal{Journal of Computational Physics}

\begin{document}

\begin{frontmatter}



\title{Hybrid Spectral Difference/Embedded Finite Volume Method \\
 for Conservation Laws}



\author{Jung J. Choi} 
\ead{jungchoi@umd.edu}

\address{The Department of Aerospace Engineering \\University of Maryland, College Park, MD 20742, USA}

\begin{abstract}
Recently, interests have been increasing towards applying the high-order methods to various engineering applications with complex geometries \cite{wang2012}. As a result, a family of discontinuous high-order methods, such as Discontinuous Galerkin (DG), Spectral Volume (SV) and Spectral Difference (SD) methods, is under active development. These methods provide high-order accurate solutions and are highly parallelizable due to the local solution reconstruction within each element. But, these methods suffer from the Gibbs phenomena when discontinuities are present in the flow fields.  Various types of limiters  \cite{du2014,zhong2013,yang2009} and artificial viscosity \cite{premasuthan2014a, persson2006} have been employed to overcome this problem. 
  
A novel hybrid spectral difference/embedded finite volume method is introduced in order to apply a discontinuous high-order method for large scale engineering applications involving discontinuities in the flows with complex geometries.  In the proposed hybrid approach, the finite volume (FV) element, consisting of structured FV subcells, is embedded in the base hexahedral element containing discontinuity, and an FV based high-order shock-capturing scheme is employed  to overcome the Gibbs phenomena. Thus, a discontinuity is captured at the resolution of FV subcells within an embedded FV element. In the smooth flow region, the SD element is used in the base hexahedral element. Then, the governing equations are solved by the SD method.  The SD method is chosen for its low numerical dissipation and computational efficiency preserving high-order accurate solutions. The coupling between the SD element and the FV element is achieved by the globally conserved mortar method \cite{kopriva1996}.  In this paper, the $5^{th}$-order WENO scheme with the characteristic decomposition is employed as the shock-capturing scheme in the embedded FV element, and the $5^{th}$-order SD method is used in the smooth flow field.

The order of accuracy study and various 1D and 2D test cases are carried out, which involve the discontinuities and vortex flows. Overall, it is shown that the proposed hybrid method results in comparable or better simulation results compared with the standalone WENO scheme when the same number of solution DOF is considered in both SD and FV elements.

\end{abstract}


\end{frontmatter}


\section{Introduction}
\label{intro}
 
The interaction of turbulence and discontinuities, e.g. shock, detonation and contact surface, in the high speed flows is commonly encountered in many engineering applications. Turbulence interacting with shock formed by a high speed aircraft generates noise downstream of the shock, which  travels to the ground and is called the sonic boom.  In such case, large scale turbulent motion interacting with shock can intensify or thicken shock, which in turn affects the level of noise downstream of the shock. Fundamental study of shock/turbulence interaction to understand the physics of the sonic boom and its reduction has been given a great attention theoretically and numerically \cite{giddings2001, lee1993}.    The shock/shear layer interaction in the high speed turbulent jet from the aircraft engine also produces screech noise, and the reduction of such noise without compensating engine performance is of great interests \cite{rona2004}.  For high speed propulsion systems, such as ramjet/scramjet, shock interacting with turbulent flows at high speed can enhance fuel/air mixing for a stable combustion with the aid of flame holding apparatus \cite{ghodke2011}.  The interaction of shock and turbulent flame plays a crucial role in the deflagration to detonation transition \cite{oran2007}.  

Turbulence is a vortex dominated chaotic process with a wide range of spatial and temporal scales while  shock is in the length scale of the molecular mean free path and is approximated by a mathematical discontinuity.  These physical processes occurring at different length scales pose  conflicting numerical requirements for successful simulations. In order to capture a wide range of length scales in turbulent flows, a numerical scheme with low numerical dissipation is required especially to capture the high wave number spectrum. On the other hand, the introduction of some numerical dissipation is needed for the numerical stability to capture discontinuities in the flow fields on a computational grid appropriate for resolving turbulent length scale.  Therefore, devising a stable and accurate numerical method for simulating turbulence interacting with discontinuities in the flow fields is a challenging task.

Various low-dissipation high-order methods have been extensively developed to reduce the numerical dissipation in the smooth or turbulent flow fields yet provide the discontinuity-capturing capability with some numerical dissipation introduced in the vicinity of discontinuity. Among various high-order schemes, weighted essentially non-oscillatory (WENO) type scheme \cite{jiang1996, shu1998}  and spectral-like compact scheme \cite{lele1992} are widely employed.   The WENO family schemes are mainly shock-capturing schemes where the order of accuracy is reduced at the location of discontinuity while the order of accuracy is preserved in the smooth flow field. However, it has been noted that the WENO scheme becomes dissipative in the smooth or turbulent flow.  The compact scheme is designed to have  very low numerical dissipation resulting in spectral-like accuracy in the smooth or turbulent flow. But the compact scheme suffers from the Gibbs phenomena when discontinuity is present in the flow field. Various efforts have been made to address these issues associated with the WENO and the compact schemes since they are introduced, especially in applying these schemes for large eddy simulation (LES) and direct numerical simulation (DNS).

With the success of early WENO scheme,  a lot of efforts have been made in order to lower the numerical dissipation in the smooth or turbulent flows. Balsara and Shu \cite{balsara2000} noted that the WENO scheme is not monotonicity preserving and, thus, introduced a monotonicity preserving weighted essentially non-oscillatory (MPWENO) scheme. They employed the monotonicity preserving bounds of Suresh and Hyunh \cite{suresh1997} and showed that the $9^{th}$ or higher order MPWENO scheme has high phase accuracy, thus suitable for compressible turbulent flows.    Martin \cite{martin2006} {\it et al.} introduced a bandwidth-optimized WENO scheme for the direct numerical simulation of turbulent compressible flows.  They proposed a set of candidate stencils to be symmetric with an additional candidate stencil. Then,  the bandwidth optimized weights for the optimal stencil  are computed by minimizing the truncation error on a given grid to maintain a small dissipation error at high wavenumbers.  The simulations of incompressible and high turbulent Mach number isotropic turbulence show good agreement compared with the simulations by the $6^{th}$-order central Pade scheme and show good high wavenumber characteristics of turbulent compressible flows. For reviews of ENO and WENO schemes and other variants, see references \cite{shu1998, henrick2005, borges2008}. 
 
The compact scheme provides spectral-like accuracy with very low numerical dissipation, but it is shown that the compact scheme suffers from the Gibbs phenomena when discontinuities are present. In order to remedy this problem,  artificial viscosity/diffusivity, filter scheme and TVD limiter  are utilized in the compact scheme for discontinuities in the flow fields.  Cook and Cabot \cite{cook2004, cook2005} introduced a high wavenumber biased spectral-like artificial shear and bulk viscosity based on the strain rate for discontinuity capturing in the compact scheme. Considering supersonic reacting flows with discontinuities,  Fiorina and Lele \cite{fiorina2006}  added artificial diffusivities to energy and species equations in addition to the artificial viscosity introduced in Cook and Cabot \cite{cook2005}, and Kawai and Lele \cite{kawai2008} extended this approach for curvilinear meshes. Mani {\it et al.} \cite{mani2009} noted that the artificial bulk viscosity in Cook and Cabot \cite{cook2005} significantly damps out the sound field while the turbulent field is not affected.  Therefore, Mani {\it et al.} proposed to replace the strain rate by the dilatation and multiply it by the Heaviside function to localize the artificial bulk viscosity near the shock and showed the improvement in sound field prediction.  Yee {\it et al.} \cite{yee1999} employed the artificial compression method (ACM) switch \cite{harten1978} as a characteristic filter to stabilize the numerical solutions and minimize the numerical dissipation near the discontinuities for the  compact scheme. By extending this filter approach, Yee and Sj{\"o}green \cite{yee2007} further developed the  high order filter scheme for multiscale Navier-Stokes and magnetohydrodynamics (MHD) systems.  The alternative approach was proposed by Cockburn and Shu \cite{cockburn1994} by using the TVD limiter and was modified by Yee \cite{yee1997} reducing spurious oscillations near the discontinuities which originate from the TVD limiter.

Another class of favored approach is the hybrid method \cite{lele2009}. The hybrid method combines a high-order low-dissipative method in the smooth or turbulent flow fields and a high-order shock-capturing scheme localizing the  numerical dissipation in the proximity of discontinuities.  In order to switch between two schemes, a discontinuity detector is devised. Adams and Shariff \cite{adams1996} introduced a hybrid method combining nonconservative compact scheme and ENO scheme. Following Adams and Shariff, Pirozzoli \cite{pirozzoli2002} used the WENO scheme in the region with discontinuities while the conservative compact scheme is used in the smooth flow region. Hill and Pullin \cite{hill2004} proposed a tuned hybrid center-difference/WENO method by developing a tuned centered difference scheme for the bandwidth optimization and coupling it with the WENO scheme for large eddy simulations with strong shocks. In these hybrid methods, it is noted that the numerical scheme to compute the convective flux is switched between the low dissipative scheme and the shock-capturing scheme depending on the characteristics of the flow. It has been shown that the hybrid scheme produces more favorable results compared with the approaches employing artificial viscosity or characteristic-based filter method \cite{lele2009, johnsen2010}.

With significant developments in the high-order methods, there have been increased interests in applying the high-order methods to engineering applications with complex geometries \cite{wang2012}. In this regard, a family of discontinuous high-order methods, such as Discontinuous Galerkin (DG) \cite{cockburn1989a,   hesthaven2002}, Spectral Volume (SV) \cite{wang2002a,   liu2006, kannan2011, kannan2012} and Spectral Difference (SD) methods \cite{liu2006a, wang2007}, has been given a great attention. Recently, Huynh \cite{huynh2007} proposed the flux reconstruction scheme, and it has been shown that the flux reconstruction scheme is a stable and efficient higher order method \cite{vincent2011, liang2013, li2013}.  The advantages of discontinuous high-order methods are that the solution reconstruction is highly localized by using a polynomial basis function within each element, thus highly parallelizable, and that these methods result in high-order accurate solutions. Due to these properties of the discontinuous high-order methods, it is well suited for the simulations using unstructured grids while applying  WENO and compact schemes may become cumbersome for unstructured grids. However, it has been shown that the family of discontinuous high-order methods also suffers from the Gibbs phenomena as in the compact scheme when discontinuities are present. To remedy this problem, limiters and artificial viscosity are employed in various studies. Du {\it et al.} \cite{du2014} and Zhong {\it et al.} \cite{zhong2013} used the WENO limiter for the correction procedure via reconstruction (CPR) method and the Runge-Kutta discontinuous Galerkin (RKDG). Yang and Wang \cite{yang2009} introduced the parameter-free generalized moment (PFGM) limiter approach for the SD method extending Cockburn and Shu's TVB limiter \cite{cockburn1994}.  Premasuthan {\it et al.} \cite{premasuthan2014a} employed the artificial viscosity approach by Cook and Cabot \cite{cook2005} for the shock-capturing in the SD method on the unstructured grid and further extended it to the adaptive mesh refinement \cite{premasuthan2014b}. 

Although the approaches employing limiters and artificial viscosity have been able to capture the discontinuities successfully in various discontinuous high-order methods, it has been noted that the discontinuities are under-resolved over a few elements while the length scale of the discontinuities may be smaller than the size of an element. In order to address this issue,  subcell shock-capturing approaches have been introduced in some of recent studies, mostly employing the DG method.  Persson and Peraire \cite{persson2006} proposed a subcell shock-capturing in  the DG method  by introducing the element-wise constant polynomial order dependent artificial viscosity by the  additional dissipative model term to the hyperbolic equation. Casoni {\it et al.} \cite{casoni2013} combined the idea of the slope limiter with the artificial viscosity approach in DG and introduced an expression to compute the artificial diffusion coefficient. Then, the piecewise constant artificial viscosity is evaluated using this expression within the shock-containing elements to capture the shock at subcell resolution. For the DG elements with non-smooth solutions, Huerta {\it et al.} \cite{huerta2012} partitioned these elements with nonoverlapping subcells and constructed the piecewise constant subcell shape functions. Then, with the piecewise constant subcell shape function, the non-smooth solutions in subcells are approximated by the first order FV shock-capturing scheme.  In Sonntag and Munz \cite{sonntag2014}, finite volume subcells are constructed in the DG elements with discontinuity keeping the Gauss solution points in the DG element as the solution points in FV subcells. Then, a Riemann solver is used for shock-capturing in FV subcells.  The flux coupling between the DG element and the FV subcells are achieved through the boundary integral in the DG formulation. Dumbser {\it et al.} \cite{dumbser2014} introduced a posteriori subcell limiting for the DG method by employing the multi-dimensional optimal order detection (MOOD) approach \cite{clain2011}. In this approach, the validity of the approximated solutions in DG elements is checked by the MOOD detection procedure. For the non-smooth solutions that require limiting, the solutions are projected onto the number of structured subcells created within this DG element and updated by the ADER-WENO method \cite{titarev2005, dumbser2013} in subcells.  Then, the updated solutions are gathered back to the DG solution points by the subcell reconstruction.

In an effort to applying the discontinuous high-order methods for large scale engineering applications involving the interaction of  turbulence and discontinuities (such as shock and detonation) in complex geometries, a novel hybrid spectral difference (SD)/embedded finite volume (FV) method is introduced in this paper for the hexahedral grids, which serve as the base grid in this method. The choice of the SD method is made based on the computational efficiency of the method compared with the DG or the SV method and its simple standard form of numerical formulations without resorting to its variants, e.g. flux reconstruction method.   In the proposed approach, the SD method is used in the smooth flow region by allocating the SD element in the base hexahedral grids.  In the non-smooth flow region where discontinuities are present, the finite volume (FV) element, consisting of structured finite volume subcells, is embedded in the base hexahedral grid in place of the SD element.  A high-order FV (or finite difference) discontinuity-capturing scheme is employed for those embedded FV elements to capture the discontinuities within an embedded FV element.  The flux coupling at the common interface between the SD element and the embedded FV element is achieved by the mortar method \cite{kopriva1996}, which satisfies the global conservation.  The advantages of the proposed hybrid method are that (a) the discontinuities are captured at the resolution of the FV subcells in the embedded FV element, resulting in sharp capturing of discontinuities within one or two  hexahedral elements depending on the nature of discontinuities (e.g. shock, detonation or contact surface) and the location of discontinuity in the hexahedral element, (b) the number of FV subcells in an FV element can be varied for high-resolution discontinuity-capturing, (c) without having to devise new schemes, current FV-based state-of-the-art high-order shock-capturing schemes can be readily used in the embedded FV element and (d) the proposed hybrid method still keeps the spatial compactness of the overall scheme  for the efficient parallel implementation since any information needed to compute numerical fluxes depends only on the immediate neighbors of the hexahedral element.  The proposed approach can be applied directly to the unstructured quadrilateral/hexahedral base grids as presented in this paper. However, in order to focus on delivering the main idea of the proposed hybrid method, the approach is described for the curvilinear structured hexahedral elements in this paper with the extended formulation of the mortar method for the hexahedral elements. Thus, for 1D and 2D test cases, the base hexahedral elements are arranged in the $x$ direction or in the $x-y$ plane, respectively, with one hexahedral element in $z$ direction.  The test results are also compared with the results obtained by Sonntag and Munz \cite{sonntag2014} and Dumbser {\it et al.} \cite{dumbser2014} as the proposed hybrid method is closely related to their works.

The paper is organized as follows. First, the system of governing equations for the conservation laws is introduced in the physical space in Section \ref{gov_eqns}.  Then, the proposed hybrid method is introduced in Section \ref{hybrid_method} with the brief descriptions of the SD method and the $5^{th}$-order WENO scheme with the characteristic decomposition and the flux coupling between the SD element and the FV element.  In Section \ref{order_accuracy}, the overall order of accuracy of the hybrid method is examined. Standard test cases in 1D and 2D are presented in Section \ref{test_cases}, and the results are discussed.    Section \ref{conclusion} summarizes the results of the current study and discusses the future works for further improvement of the proposed hybrid method.

\section{Governing Equations}
\label{gov_eqns}
A system of 3D unsteady compressible Euler equation is considered and written in the conservation form as,

\begin{equation}
\frac{\partial \bf{Q}}{\partial t} + \frac{\partial \bf{F}}{\partial x}  + \frac{\partial \bf{G}}{\partial y}  + \frac{\partial \bf{H}}{\partial z} =0, \label{phys_euler}
\end{equation} 
where $ \mbox{\bf Q} = [\rho,~\rho u,~\rho v,~\rho w,~\rho E]^T$. The inviscid flux vectors are
\begin{equation}
\bf{F} = \left(
\begin{array}{l}
\rho u \\
\rho uu + p \\
\rho uv \\
\rho uw \\
\rho u (E+p/ \rho)
\end{array} \right),  \qquad
\bf{G} = \left(
\begin{array}{l}
\rho v \\
\rho vu \\
\rho vv + p \\
\rho vw \\
\rho v (E+p/ \rho)
\end{array} \right),  \qquad
\bf{H} = \left(
\begin{array}{l}
\rho w \\
\rho wu \\
\rho wv \\
\rho ww + p\\
\rho w (E+p/ \rho)
\end{array} \right), 
\end{equation}
where $\rho$ is the density, $p$ is the pressure, and $u$, $v$ and $w$ are the velocity components. The total energy is denoted by $E$ and defined as 
\begin{equation}
E = e + \frac{1}{2} (u^2+v^2+w^2),
\end{equation} 
where $e$ is the internal energy, and $e=C_v T$ under the calorically perfect gas assumption. The system of governing equations is closed with the ideal gas equation, $p =\rho R T$, where $R$ is the specific gas constant and $T$ is the temperature.

\begin{figure}[tb]
\begin{center}
 \includegraphics[width=6in]{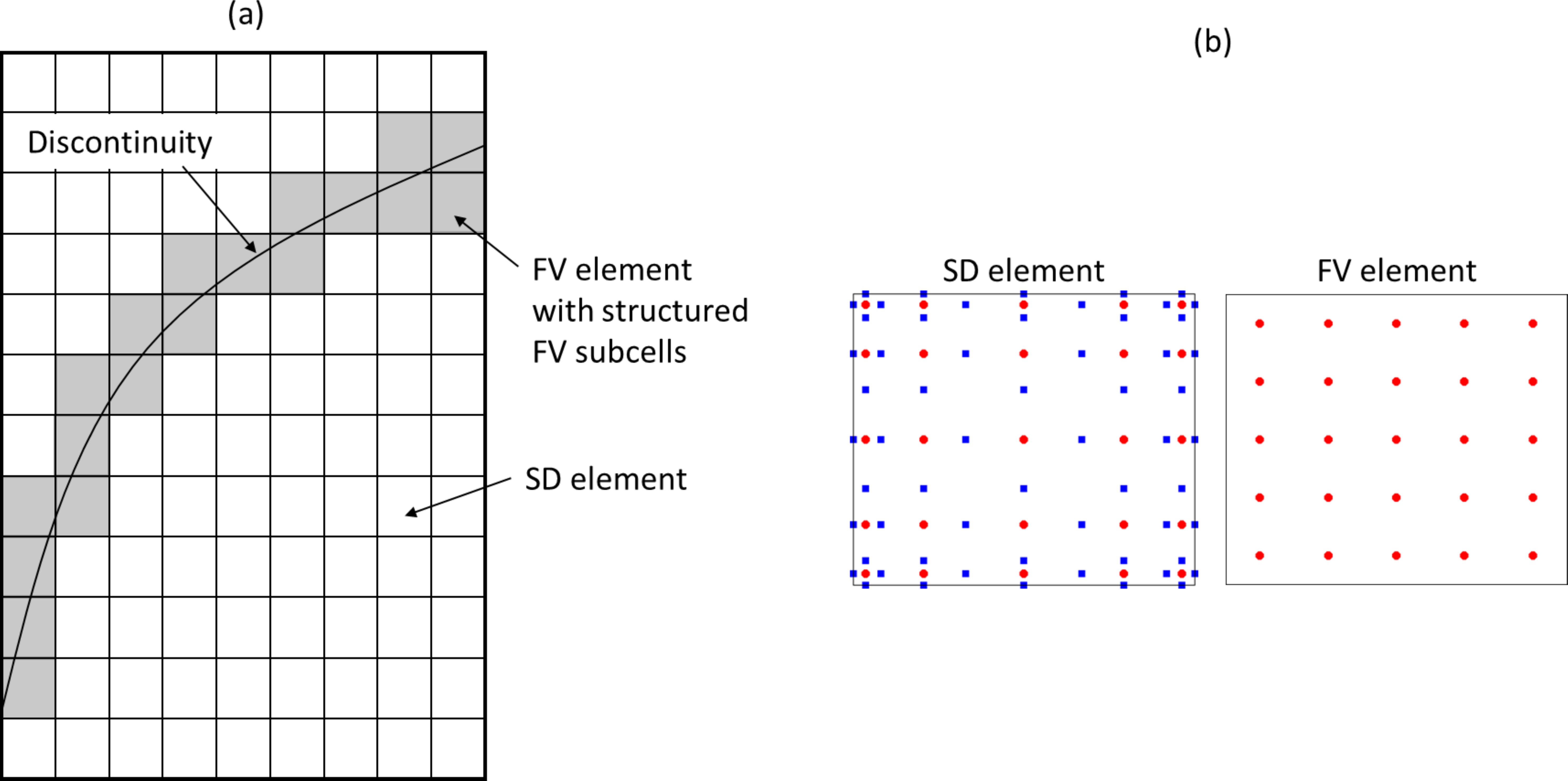}  
 \caption{Schematic diagram of embedded FV element in the hybrid method. (a) Overall schematic diagram where FV element is embedded in the base elements containing discontinuity (depicted by gray colored elements)  and (b) the locations of quadrature points for solution (red filled circle) and flux (blue filled square) in the SD element and the cell centered subcell solution points (red filled circle) in the FV element.}
 \label{sdweno}
 \end{center}
\end{figure}

\section{Hybrid Method}
\label{hybrid_method}

The main idea of the proposed hybrid approach is to utilize the finite volume based high-order shock-capturing scheme locally in the region where the discontinuity is present by dynamically embedding finite volume subcells within the base hexahedral element, as depicted in Fig. \ref{sdweno}(a). Away from the discontinuity where the flow field is smooth, the SD element is used in the base hexahedral element where the SD method is employed to solve the governing equations.   As mentioned in Section \ref{intro}, the SD method is chosen for its computational efficiency and simple numerical formulations.   For the finite volume based high-order shock-capturing method, the standard WENO scheme with the characteristic decomposition is employed while other high-order shock-capturing schemes can also be used.  For the purpose of presenting main idea of the proposed hybrid method, the $5^{th}$-order SD method and the $5^{th}$-order WENO scheme with characteristic decomposition are employed in the present study. The same order of accuracy is chosen such that the overall order of the hybrid method in space is to be the 5$^{th}$ order. Choosing the number of FV subcells in the embedded FV element, as discussed Dumbser {\it et al.} \cite{dumbser2014}, is important for determining  the optimal time step for both SD method and high-order shock-capturing scheme on the FV subcells in the embedded FV element and reducing the truncation error of the high-order scheme on the FV subcells.  It is suggested in \cite{dumbser2014} that the optimal number of FV subcells is $2N+1$ per spatial dimension, where $N$ is the order of polynomial used in the DG element (or the SD element in this paper), but it is also mentioned that the choice of the number of FV subcells is left to the user.  In the current study, the number of FV subcells used in the FV element is $5^3$. This choice is suboptimal, but it is intended in order to have the same solution degree of freedom (DOF) in the FV element as in the $5^{th}$-order SD element and to investigate the effect of the suboptimal number of FV subcells on the solution in the hybrid method.  The location of solution points is shown in Fig. \ref{sdweno}(b) for the SD element and the FV element with structured FV subcells.  The essential part of the current hybrid method lies in computing the common interface flux coupling between the SD element and the FV element, which is achieved by the mortar method \cite{kopriva1996}.  

Non-uniform hexahedral element in physical space is transformed to the standard unit cubic element in the computational domain by the geometric mapping.   Accordingly, the system of governing equations,  Eq. (\ref{phys_euler}), in the physical domain is transformed to the computational domain using the Jacobian transformation matrix,

\begin{equation}
\frac{\partial \tilde{\bf{Q}}}{\partial t} + \frac{\partial \tilde{\bf{F}}}{\partial \xi}  + \frac{\partial \tilde{\bf{G}}}{\partial \eta}  + \frac{\partial \tilde{\bf{H}}}{\partial \zeta} =0, \label{comp_euler}
\end{equation} 
where $\tilde{\bf{Q}} = |\bf{J}| \bf{Q}$ with the Jacobian matrix, {\bf J}, and
\begin{equation}
\left( \begin{array}{l}
\tilde{\bf{F}} \\
\tilde{\bf{G}} \\
\tilde{\bf{H}} 
\end{array} \right)
= |\bf{J}| J^{-1}
\left( \begin{array}{l}
\bf{F} \\
\bf{G} \\
\bf{H} 
\end{array} \right).
\end{equation}

The time integration in both SD and WENO methods is carried out by the $3^{rd}$-order strong stability preserving (SSP) Runge-Kutta (RK) scheme with low storage method \cite{spiteri2002}, and it is given in the general form as
\begin{eqnarray}
\tilde{\bf{Q}}^{(0)} &=& \tilde{\bf{Q}}^{n}, \qquad d\tilde{\bf{Q}}^{(0)} = 0,  \nonumber  \\
d\tilde{\bf{Q}}^{(i)} &=& A_i d\tilde{\bf{Q}}^{(i-1)} + \Delta t L(\tilde{\bf{Q}}^{(i-1)}), \qquad i = 1, ..., m, \\
\tilde{\bf{Q}}^{(i)} &=& \tilde{\bf{Q}}^{(i-1)} + B_i d\tilde{\bf{Q}}^{(i-1)},  \qquad \qquad ~~~~i= 1,...,m,  \nonumber \\
\tilde{\bf{Q}}^{n+1} &=& \tilde{\bf{Q}}^{(m)},  \nonumber
\end{eqnarray}
with  $A_1=0$, and $m=3$ for 3 stage. The coefficients $A_i$ and $B_i$ are given in \cite{spiteri2002}. 
 
In the following, the brief description of the standard SD method and the $5^{th}$-order WENO scheme with characteristic decomposition is given followed by the reconstruction of the common interface flux between the SD element and the FV element by the mortar method.

\subsection{Spectral difference method}
\label{sd_method}
In the spectral difference method, two sets of points are required inside the SD element, namely solution points and flux points. The unknown solutions are located at the solution points while flux values are located at the flux points. The solution points are defined by the Gauss points given by
\begin{equation}
X_s = \frac{1}{2} \left[ 1 - \cos \left( \frac{2s-1}{2n} \pi \right) \right], \quad s = 1,2, \dots, n.  \label{gauss_pts}
\end{equation}

For the flux points, the Legendre-Gauss points with additional two end points, 0 and 1, are used. The Legendre-Gauss polynomial is given by,
\begin{equation}
P_s(\xi) = \frac{2s-1}{s} \left(2 \xi -1 \right) P_{s-1}(\xi) - \frac{s-1}{s}P_{s-2}(\xi), \quad s = 1, ..., n-1,
\label{legendre_gauss_pts}
\end{equation}
with $P_{-1}(\xi)=0$ and $P_{0}(\xi)=1$. The flux points are then defined with the roots of the Legendre-Gauss polynomial $P_{n-1}$ with two end points for the $n^{th}$-order SD element.

The solution points are used to construct $(n-1)$ degree Lagrange polynomial given by
\begin{equation}
h_i(X) = \prod_{s=1,s \ne i}^{n} \left( \frac{X-X_s}{X_i-X_s} \right), \label{hpoly}
\end{equation}
for solution reconstruction. Likewise, $n$ degree Lagrange polynomial is constructed for the flux reconstruction given by
\begin{equation}
l_{i+1/2}(X) = \prod_{s=0,s \ne i}^{n} \left( \frac{X-X_{s+1/2}}{X_{i+1/2} - X_{s+1/2}} \right). \label{lpoly}
\end{equation}

Then, the solution is reconstructed by the tensor product of one-dimensional Lagrange polynomials. 
\begin{equation}
{\bf Q}(\xi,\eta,\zeta) = \sum_{k=1}^{n} \sum_{j=1}^{n} \sum_{i=1}^{n} \frac{\tilde{\bf{Q}}_{i,j,k}}{|{\bf J}_{i,j,k}|} h_i(\xi) \cdot h_j(\eta) \cdot h_k(\zeta).  \label{q_interp}
\end{equation}
And the flux reconstruction is given by
\begin{eqnarray}
\tilde{\bf{F}}(\xi,\eta,\zeta) &=& \sum_{k=0}^{n} \sum_{j=0}^{n} \sum_{i=0}^{n} \tilde{\bf{F}}_{i+1/2,j,k} l_{i+1/2}(\xi) \cdot h_j(\eta) \cdot h_k(\zeta), \label{f_interp} \\
\tilde{\bf{G}}(\xi,\eta,\zeta) &=& \sum_{k=0}^{n} \sum_{j=0}^{n} \sum_{i=0}^{n} \tilde{\bf{G}}_{i+1/2,j,k} h_i(\xi) \cdot l_{j+1/2}(\eta) \cdot  h_k(\zeta), \label{g_interp} \\
\tilde{\bf{H}}(\xi,\eta,\zeta) &=& \sum_{k=0}^{n} \sum_{j=0}^{n} \sum_{i=0}^{n} \tilde{\bf{H}}_{i+1/2,j,k}  h_i(\xi) \cdot h_j(\eta) \cdot l_{k+1/2}(\zeta). \label{h_interp}
\end{eqnarray}

In advancing the solution in time, the conservative variable, $\tilde{\bf Q}$, is interpolated to the flux points using Eq. (\ref{q_interp}). The numerical fluxes are constructed at the flux points using Eqs. (\ref{f_interp})-(\ref{h_interp}) with $\tilde{\bf Q}$ at flux points. Note that the numerical fluxes reconstructed are only element-wise continuous, but discontinuous between neighboring elements. Thus, it requires a Riemann solver, such as the Rusanov solver \cite{rusanov1961} or the AUSM$^+$-up scheme \cite{liou2006ausm}, to compute the common fluxes at the element interface.  Then, the derivatives of the fluxes at the solution points are computed by
\begin{eqnarray}
\left( \frac{\partial \tilde{\bf F}}{\partial \xi} \right)_{i,j,k} &=& \sum_{m=0}^{n} \tilde{\bf F}_{m+1/2,j,k} \cdot l'_{m+1/2}(\xi_i), \label{f_deriv}\\
\left( \frac{\partial \tilde{\bf G}}{\partial \xi} \right)_{i,j,k} &=& \sum_{m=0}^{n} \tilde{\bf G}_{i,m+1/2,k} \cdot l'_{m+1/2}(\eta_j), \label{g_deriv}\\
\left( \frac{\partial \tilde{\bf H}}{\partial \xi} \right)_{i,j,k} &=& \sum_{m=0}^{n} \tilde{\bf H}_{i,j,m+1/2} \cdot l'_{m+1/2}(\zeta_k), \label{h_deriv}
\end{eqnarray}
where $l'$ is the spatial derivative of the Lagrange polynomial.  And finally, the solution is updated by
\begin{equation}
\frac{\partial \tilde{\bf Q}_{i,j,k}}{\partial t} = - \left( \frac{\partial \tilde{\bf F}}{\partial \xi} + \frac{\partial \tilde{\bf G}}{\partial \eta} +  \frac{\partial \tilde{\bf H}}{\partial \zeta} \right)_{i,j,k} \label{q_update}
\end{equation}
 at each solution point.

\begin{figure}[tb]
\begin{center}
 \includegraphics[width=3.2in]{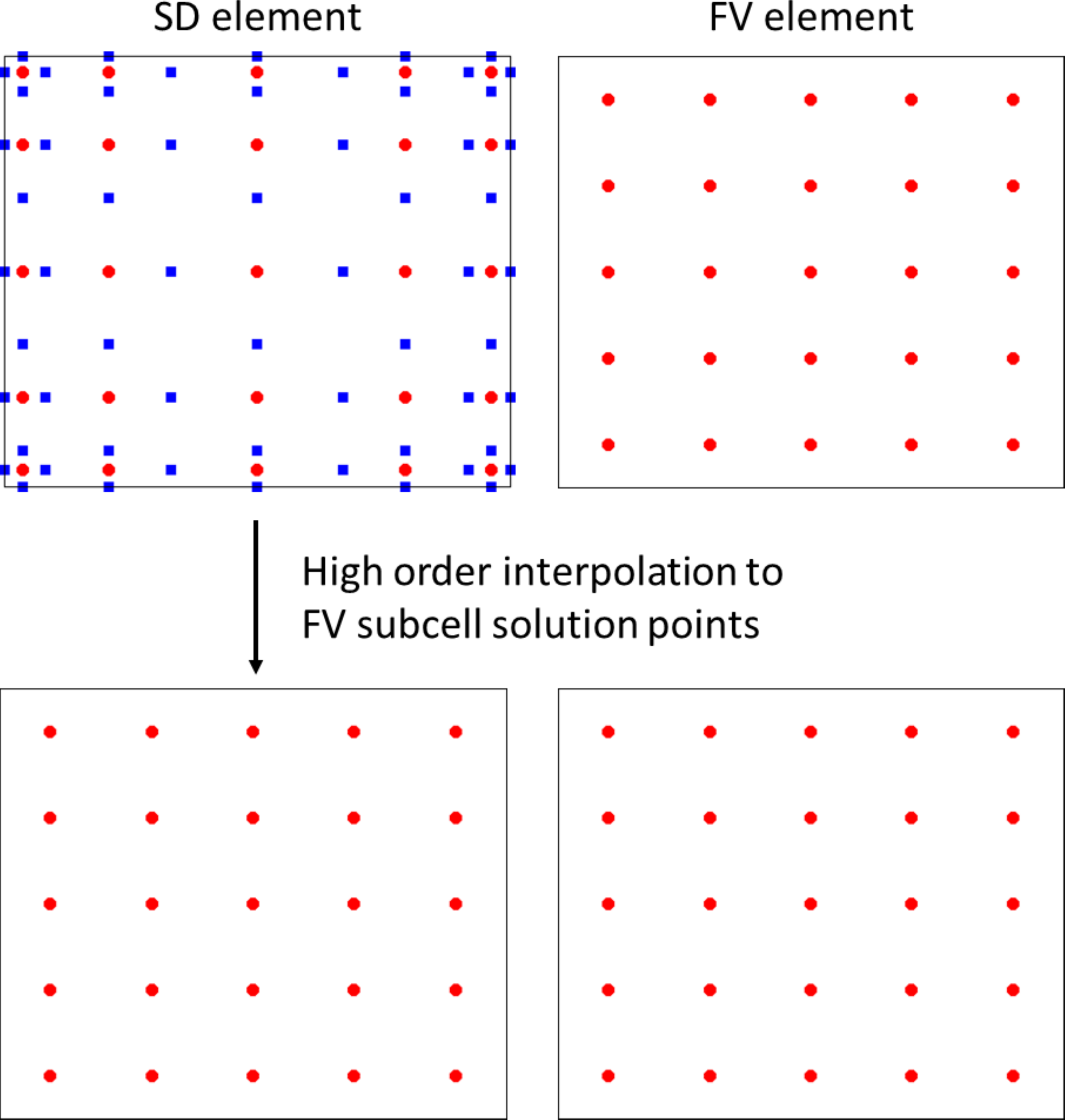}
 \caption{Interpolation of conservative variables to FV subcell solution points in the neighboring SD element of the  FV element, which are used in the WENO interpolation.}
 \label{weno_interp}
 \end{center}
\end{figure}

\begin{figure}[t]
\begin{center}
 \includegraphics[width=4in]{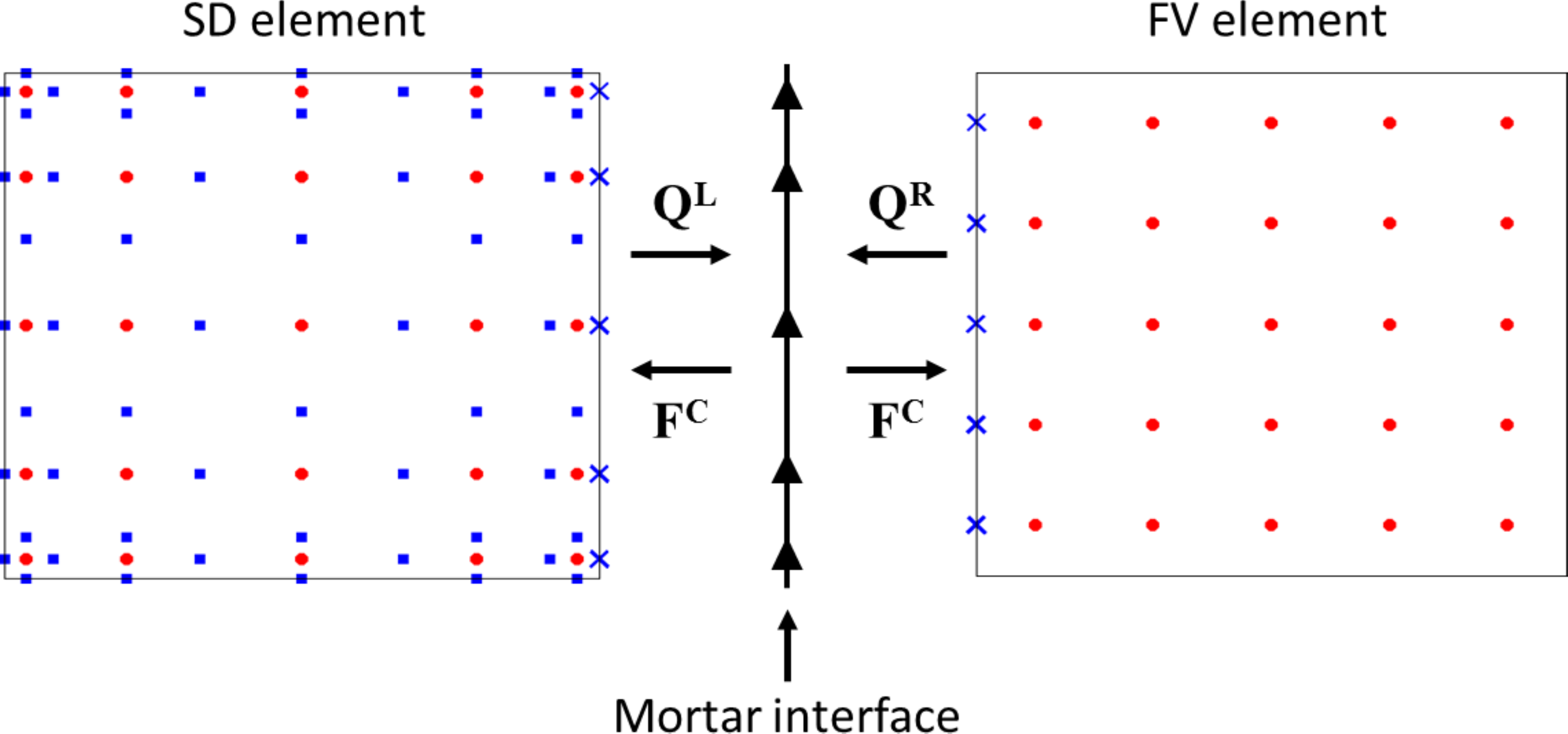}
 \caption{Computation of common interface flux between the SD element and the FV element using the mortar interface. }
 \label{coupling}
 \end{center}
\end{figure}

\subsection{$5^{th}$-order WENO scheme}
\label{weno_scheme}
In the WENO scheme, the numerical flux at the mid point, $\tilde{\bf F}_{i+1/2,j,k}$ is constructed from the left and right states by the weighted averaging procedure and using a Riemann solver as following,
\begin{equation}
\tilde{\bf F}_{i+1/2,j,k} = {\cal R} \left(\tilde{\bf Q}_{i+1/2,j,k}^{L},\tilde{\bf Q}_{i+1/2,j,k}^{R} \right), 
\end{equation}
where $\tilde{\bf Q}_{i+1/2,j,k}^{L}$ and $\tilde{\bf Q}_{i+1/2,j,k}^{R}$ are obtained from the WENO interpolation, and $\cal R$ is an operator denoting a Riemann solver. For brevity, only the WENO interpolation for the left state in $i$ direction is described. The right state can be computed in the same way using the symmetric arrangement of the stencils. The conservative variable, $\tilde{\bf Q}$,  is transformed into the characteristic variable as follows
\begin{equation}
{\cal Q}_{k,m}^{L} = l_{i+1/2,m} \tilde{\bf Q}_{k}^{L} \quad (k=i-2,i-1,i,i+1,i+2),
\end{equation}
where ${ \bf \cal Q}_{k,m}$ denotes the $m$th characteristic variable, and $l_{i+1/2,m}$ denotes the $m^{th}$ left eigenvector of the matrix $\partial \tilde{\bf F} / \partial \tilde{\bf Q}$ at $i+1/2$, which is computed with the Roe average of values at $i$ and $i+1$. Then, the characteristic variable, ${\bf \cal Q}_{i+1/2,m}^{L}$, is constructed as follows
\begin{equation}
{\bf \cal Q}_{i+1/2,m}^{L} = w_{i,m}^{1} {\bf \cal Q}_{i+1/2,m}^{L,1} + w_{i,m}^{2} {\bf \cal Q}_{i+1/2,m}^{L,2} + w_{i,m}^{3} {\bf \cal Q}_{i+1/2,m}^{L,3},
\end{equation}  
where 
\begin{eqnarray}
{ \bf \cal Q}_{i+1/2,m}^{L,1} &=& \frac{1}{3} {\bf \cal Q}_{i-2,m}^{L} - \frac{7}{6} {\bf \cal Q}_{i-1,m}^{L} +\frac{11}{6} {\bf \cal Q}_{i,m}^{L}, \\
{ \bf \cal Q}_{i+1/2,m}^{L,2} &=& -\frac{1}{6} {\bf \cal Q}_{i-1,m}^{L} + \frac{5}{6} {\bf \cal Q}_{i,m}^{L} +\frac{1}{3} {\bf \cal Q}_{i+1,m}^{L}, \\
{ \bf \cal Q}_{i+1/2,m}^{L,3} &=& \frac{1}{3} {\bf \cal Q}_{i,m}^{L} + \frac{5}{6} {\bf \cal Q}_{i+1,m}^{L} -\frac{1}{6} {\bf \cal Q}_{i+2,m}^{L}, 
\end{eqnarray}
and $w^1$, $w^2$ and $w^3$ are the nonlinear weights, which are determined as
\begin{equation}
w_{i,m}^{k} = \frac{\alpha_{i,m}^{k}}{\alpha_{i,m}^{1} + \alpha_{i,m}^{2} + \alpha_{i,m}^{3}}, \quad k=1,2,3,
\end{equation}
where 
\begin{equation}
\alpha_{i,m}^k = \frac{C^k}{(IS_{i,m}^k + \epsilon)^p} \quad k=1,2,3,
\end{equation}
with 
\begin{equation}
C^1 = \frac{1}{10}, \qquad C^2 = \frac{6}{10}, \qquad C^3 = \frac{3}{10}, 
\end{equation}
and
\begin{eqnarray}
IS_{i,m}^{1} &=& \frac{1}{4} \left( {\bf \cal Q}_{i-2,m}^{L} - 4 {\bf \cal Q}_{i-1,m}^{L} + 3 {\bf \cal Q}_{i,m}^{L} \right)^2  + \frac{13}{12} \left( {\bf \cal Q}_{i-2,m}^{L} - 2 {\bf \cal Q}_{i-1,m}^{L} +  {\bf \cal Q}_{i,m}^{L} \right)^2, \\
IS_{i,m}^{2} &=& \frac{1}{4} \left( -{\bf \cal Q}_{i-1,m}^{L} + {\bf \cal Q}_{i+1,m}^{L} \right)^2  + \frac{13}{12} \left( {\bf \cal Q}_{i-1,m}^{L} - 2 {\bf \cal Q}_{i,m}^{L} +  {\bf \cal Q}_{i+1,m}^{L} \right)^2,, \\
IS_{i,m}^{3} &=& \frac{1}{4} \left( -3{\bf \cal Q}_{i,m}^{L} + 4 {\bf \cal Q}_{i+1,m}^{L} - {\bf \cal Q}_{i+2,m}^{L} \right)^2  + \frac{13}{12} \left( {\bf \cal Q}_{i,m}^{L} - 2 {\bf \cal Q}_{i+1,m}^{L} +  {\bf \cal Q}_{i+2,m}^{L} \right)^2.
\end{eqnarray}
Here, $\epsilon=10^{-6}$ and $p=2$. Then, the characteristic variable,  ${\bf \cal Q}_{i+1/2,m}^{L}$, is transformed back to the conservative variable, $\tilde{\bf Q}$, as follows
\begin{equation}
\tilde{\bf Q}_{i+1/2}^{L} = \sum_m {\bf \cal Q}_{i+1/2,m}^{L}r_{i+1/2,m},
\end{equation}
where $r_{i+1/2,m}$ is the $m^{th}$ right eigenvector of the matrix $\partial \tilde{\bf F} / \partial \tilde{\bf Q}$ at $i+1/2$.
 
Then, the derivative of flux is evaluated as,
\begin{equation}
\left( \frac{\partial \tilde{\bf F}}{\partial \xi} \right)_{i,j,k} = \frac{1}{\Delta \xi} \left( \tilde{\bf F}_{i+1/2,j,k} -\tilde{\bf F}_{i-1/2,j,k} \right).
\end{equation} 
Likewise, $\left( \frac{\partial \tilde{\bf G}}{\partial \eta} \right)_{i,j,k}$ and $\left( \frac{\partial \tilde{\bf H}}{\partial \zeta} \right)_{i,j,k}$ are computed. Then, the solution is advanced using Eq. (\ref{q_update}).

In the case that the neighboring element is the SD element as shown in Fig. \ref{weno_interp}, the solutions at the Gauss points are interpolated to the solution points of FV subcells in the neighboring SD element using Eq. (\ref{q_interp}) with the same order of Lagrange polynomial. Then, those interpolated values at the FV subcell solution points in the neighboring SD element are used to construct the numerical fluxes by the WENO interpolation procedure in the FV  subcells.  It should be noted that there is the truncation error associated with this interpolation step between the Gauss points and solution points of FV subcells in the neighboring SD element. As discussed in Dumbser {\it et al.} \cite{dumbser2014}, the optimal number of FV subcells can be used to minimize the truncation error associated with this interpolation step.     
 
In some cases, such as the presence of strong discontinuities, a limiter is required in the WENO scheme for the numerical stability \cite{balsara2000}.  The following form of total variation diminishing (TVD) limiter  is employed as described in \cite{houim2011}.  
\begin{eqnarray}
\phi = \max \left[ 0, \min \left( \alpha, \alpha \frac{Q_{i+1}-Q_i}{Q_i-Q_{i-1}}, 2 \frac{\hat{Q}_{i+1/2}-Q_i}{Q_i-Q_{i-1}} \right) \right],  \label{tvd_lim}
\end{eqnarray}  
where $\hat{Q}_{i+1/2}$ is the interpolated value using the $5^{th}$-order WENO interpolation, and $\alpha$ is a constant, which is set to 2. Then, the slope limited interpolation is defined by
\begin{equation}
Q_{i+1/2} = Q_i + 0.5 \left( Q_i - Q_{i-1} \right) \phi. \label{slope_tvd}
\end{equation}

It is noted that the main purpose of using TVD limiter with the WENO scheme is to access its effect on the numerical solutions in the proposed hybrid method while the MPWENO \cite{balsara2000} may be used in place of the TVD limiter.

\subsection{SD/embedded FV interface flux reconstruction}
\label{mortar_coupling}
When the neighboring elements are the same type at the common element interface, the standard SD or the WENO interpolation procedure is used to compute the common flux at the element interface. However, as the FV elements are embedded in the region where discontinuities are present, there are element interfaces that have different types of elements in the neighbor with different sets of solution points in each element.   For these interfaces, consistent and conservative common flux should be evaluated. In this study, reconstruction of common flux at this type of interface is achieved by employing the mortar method \cite{kopriva1996}.  The mortar method is designed to satisfy the global conservation and the outflow conditions for a hyperbolic problem, where waves should pass through the element interface.  It is discussed in Kopriva \cite{kopriva1996} that the least square projection satisfies these two conditions.  The mortar projection procedure to construct common flux at the interface between hexahedral elements is developed by extending the formulation described in \cite{kopriva1996} for the interface between quadrilateral elements.  The schematic diagram of mortar procedure is shown in Fig. \ref{coupling}.

\begin{figure}[t]
\begin{center}
 \includegraphics[width=6in]{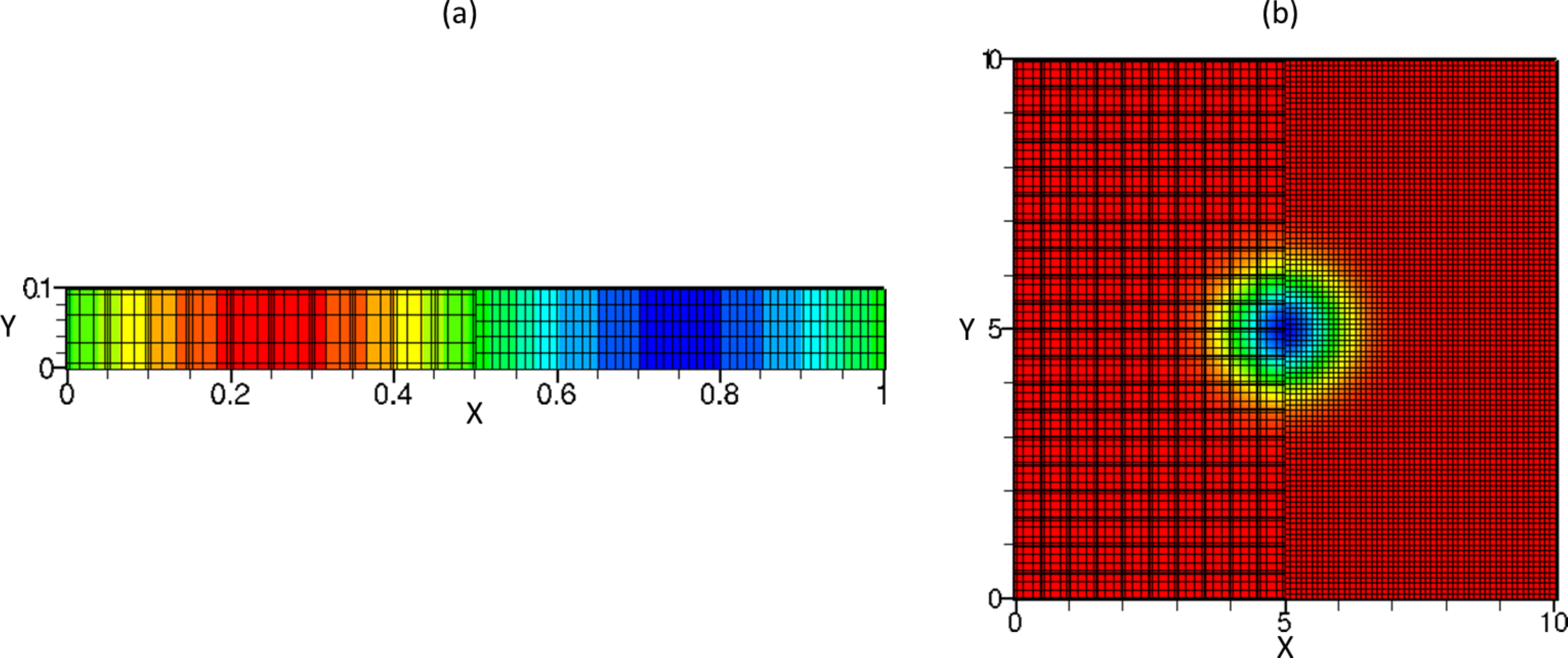} 
 \caption{Computational mesh showing static SD elements and embedded FV subcells in the FV elements with the initial density profile for (a) translating density sine wave and (b) translating vortex.}
 \label{sine_vort_grid}
 \end{center}
\end{figure}

First, the conservative variable, $\tilde{\bf Q}$, is computed at the flux points on the face of the participating element.  Referring to Fig. \ref{coupling}, the values of $\tilde{\bf Q}$  in the SD element are computed at the flux points on the right side face of the SD element, shown as blue cross mark in Fig. \ref{coupling}. In the FV elements, $\tilde{\bf Q}$ at flux points on the left side face of the FV element, denoted by blue cross mark in Fig. \ref{coupling}, is computed following the WENO procedure. A mortar interface is defined by the Gauss points with the maximum number of solution points in each spatial direction of participating face between the SD and the FV elements. Then, the Lagrange polynomial is constructed with the Gauss points on the mortar interface with the polynomial order, $J-1$, where $J= \max(M_{sd}, M_{fv})$. Here, $M_{sd}-1$ and $M_{fv}-1$ are the order of polynomials in the SD element and the FV element. The $\tilde{\bf Q}$ values from SD element and the FV element are  projected onto the mortar interface by the least square projection given by
\begin{equation}
\sum_{i=0}^{J-1} \sum_{j=0}^{J-1} \Phi_{i,j}  M_{ij, mn} = \sum_{i=0}^{M_{sd/fv}-1} \sum_{j=0}^{M_{sd/fv}-1} U_{i,j} S_{ij, mn}, \qquad m, n = 0, ...,J-1,
\end{equation}
where $\Phi_{i,j}$ is the interpolated values on the mortar interface, and $U_{i,j}$ is the $\tilde{\bf Q}$ variable on the participating face of either the SD element or the FV element.  The projection matrices, $M_{ij,mn}$ and $S_{ij,mn}$ are defined as
\begin{eqnarray}
M_{ij, mn} &=& \int_{0}^{1} \int_{0}^{1} h_{i}^{M}(\xi) h_{j}^{M}(\eta) h_{m}^{M}(\xi) h_{n}^{M}(\eta) d \xi d \eta, \\
S_{ij, mn} &=& \int_{0}^{1} \int_{0}^{1} h_{i}(\xi) h_{j}(\eta) h_{m}^{M}(\xi) h_{n}^{M}(\eta) d \xi d \eta,
\end{eqnarray}
where $h_i$ and $h_j$ are the Lagrange polynomials defined on the face of SD element and the FV element, and $h_i^M$ and $h_j^M$ are the Lagrange polynomials defined on the mortar interface. The integrals in the projection matrices, $M_{ij,mn}$ and $S_{ij,mn}$, can be evaluated by using  a Clenshaw$-$Curtis quadrature \cite{canuto1987}.

These projected $\tilde{\bf Q}$ values to the mortar interface from both SD and FV elements then serve as left and right states for the Riemann solver.  The common interface flux is computed by the Riemann solver at the mortar flux points.  Finally, the common flux at the mortar interface is projected back to the participating faces of the SD and FV elements by  
\begin{equation}
\sum_{i=0}^{M_{sd/fv}-1} \sum_{j=0}^{M_{sd/fv}-1} \Phi_{i,j}  M_{ij, mn} = \sum_{i=0}^{J-1} \sum_{j=0}^{J-1} U_{i,j} S_{ij, mn}, \qquad m, n = 0, ..., M_{sd/fv}-1,
\end{equation}
where $\Phi_{i,j}$ is the interpolated flux on participating face of the SD or the FV element, and $U_{i,j}$ is the flux on the mortar interface.  The projection matrices, $M_{ij,mn}$ and $S_{ij,mn}$ are defined as
\begin{eqnarray}
M_{ij, mn} &=& \int_{0}^{1} \int_{0}^{1} h_{i}(\xi) h_{j}(\eta) h_{m}(\xi) h_{n}(\eta) d \xi d \eta, \\
S_{ij, mn} &=& \int_{0}^{1} \int_{0}^{1} h_{i}^{M}(\xi) h_{j}^{M}(\eta) h_{m}(\xi) h_{n}(\eta) d \xi d \eta.
\end{eqnarray}
The projected flux at the face of the SD or the FV element is used in computing the derivative of the numerical flux. 

While the mortar method is developed for coupling two neighboring spectral elements in Kopriva \cite{kopriva1996}, it should be noted that the same least square projection formulation is applied in this work in an attempt to couple two different types of elements, namely SD element and the FV element. Therefore, it should be noted that there is a difference in the magnitude of truncation error between the SD element/mortar interface projection step and FV element/mortar interface projection step.  The truncation error in the FV element/mortar interface projection step tends to be larger than the SD element/mortar interface projection step and may affect the order of accuracy depending on the number of FV subcells used in the FV element. Using the optimal number of FV subcells discussed in Section \ref{weno_scheme} (see discussion for Fig. \ref{weno_interp}), the truncation error may be minimized in the FV element/mortar interface projection step.  The effect of truncation error associated with the mortar projection step is discussed in the order of accuracy study in Section \ref{order_accuracy}.

\subsection{Discontinuity detector}
In order to identify the element that contains the discontinuity in the flow field, the detector described in \cite{jameson1981} is employed and is given by 
\begin{equation}
\left| \frac{\phi_{i+1} - 2 \phi_i + \phi_{i-1}}{\phi_{i+1} + 2 \phi_i + \phi_{i-1}}  \right| > \epsilon_s,
\end{equation}
where $\phi$ is the pressure or the density, and $\epsilon_s$ is the user defined constant.  The value of $\epsilon_s$ is varied to keep the embedded FV regions close to the discontinuity.  
 
\begin{table}[t]
\begin{center}
\caption{Translating density sine wave: error and order of accuracy.}
\begin{tabular} { c|c|c c|c c|c c }
\hline \hline
            & Mesh                  &$L_1$ error  & Order& $L_2$ error& Order& $L_{\infty}$ error &  Order \\ \hline 
\multirow{3}{*}
{SD\_RUS }  & 20$\times$1$\times$1  & 4.3287e-08  & -    & 4.9659e-08 & -    & 1.0150e-07       & -  \\ 
            & 40$\times$1$\times$1  & 1.6219e-09  & 4.74 & 1.8590e-09 & 4.74 & 3.7739e-09       & 4.75 \\   
 		    & 80$\times$1$\times$1  & 5.1751e-11  & 4.97 & 5.9458e-11 & 4.97 & 1.0771e-10       & 5.13 \\ \hline 
\multirow{3}{*}
{SD\_AUSM}  & 20$\times$1$\times$1  & 1.2598e-08  & -    & 1.8073e-08 &  -   & 4.6771e-08       &  - \\ 
            & 40$\times$1$\times$1  & 3.9029e-10  & 5.01 & 5.6199e-10 & 5.01 & 1.4511e-09       & 5.01  \\   
 			& 80$\times$1$\times$1  & 1.5775e-11  & 4.63 & 1.9871e-11 & 4.82 & 4.6229e-11       & 4.97  \\ \hline 
\multirow{3}{*}
{WENO\_AUSM}& 100$\times$1$\times$1  & 2.1569e-07  & -    & 2.4392e-07 & -    &  4.1328e-07     & - \\ 
            & 200$\times$1$\times$1  & 6.7398e-09  & 5.00 & 7.5514e-09 & 5.01 &  1.2911e-08     & 5.00 \\   
 		    & 400$\times$1$\times$1  & 2.1319e-10  & 4.98 & 2.3769e-10 & 4.99 &  3.9799e-10     & 5.02 \\ \hline  
\multirow{4}{*}
{WENO\_AUSM\_TVD} & 100$\times$1$\times$1  & 6.5065e-04  & -    & 1.2092e-03 & -     & 3.5669e-03 & - \\ 
                  & 200$\times$1$\times$1  & 1.2804e-04  & 2.35 & 3.1984e-04 & 1.92  & 1.2412e-03 & 1.52  \\   
 				  & 400$\times$1$\times$1  & 2.5361e-05  & 2.34 & 8.2726e-05 & 1.95  & 4.2338e-04 & 1.55 \\  \hline 
\multirow{3}{*}
{HYBRID\_STATIC } 
	        & 20$\times$1$\times$1  & 4.2916e-07 & -    & 8.0529e-07 & -     & 3.5788e-06 & - \\ 
            & 40$\times$1$\times$1  & 5.4189e-08 & 2.99 & 1.2892e-07 & 2.64  & 7.7709e-07 & 2.20 \\   
 			& 80$\times$1$\times$1  & 7.4116e-09 & 2.87 & 2.1448e-08 & 2.59  & 1.6716e-07 & 2.22 \\    
\hline \hline
\end{tabular}
\label{sine_error} 
\end{center}
\end{table}

\begin{figure}[t]
\begin{center}
 \includegraphics[width=6in]{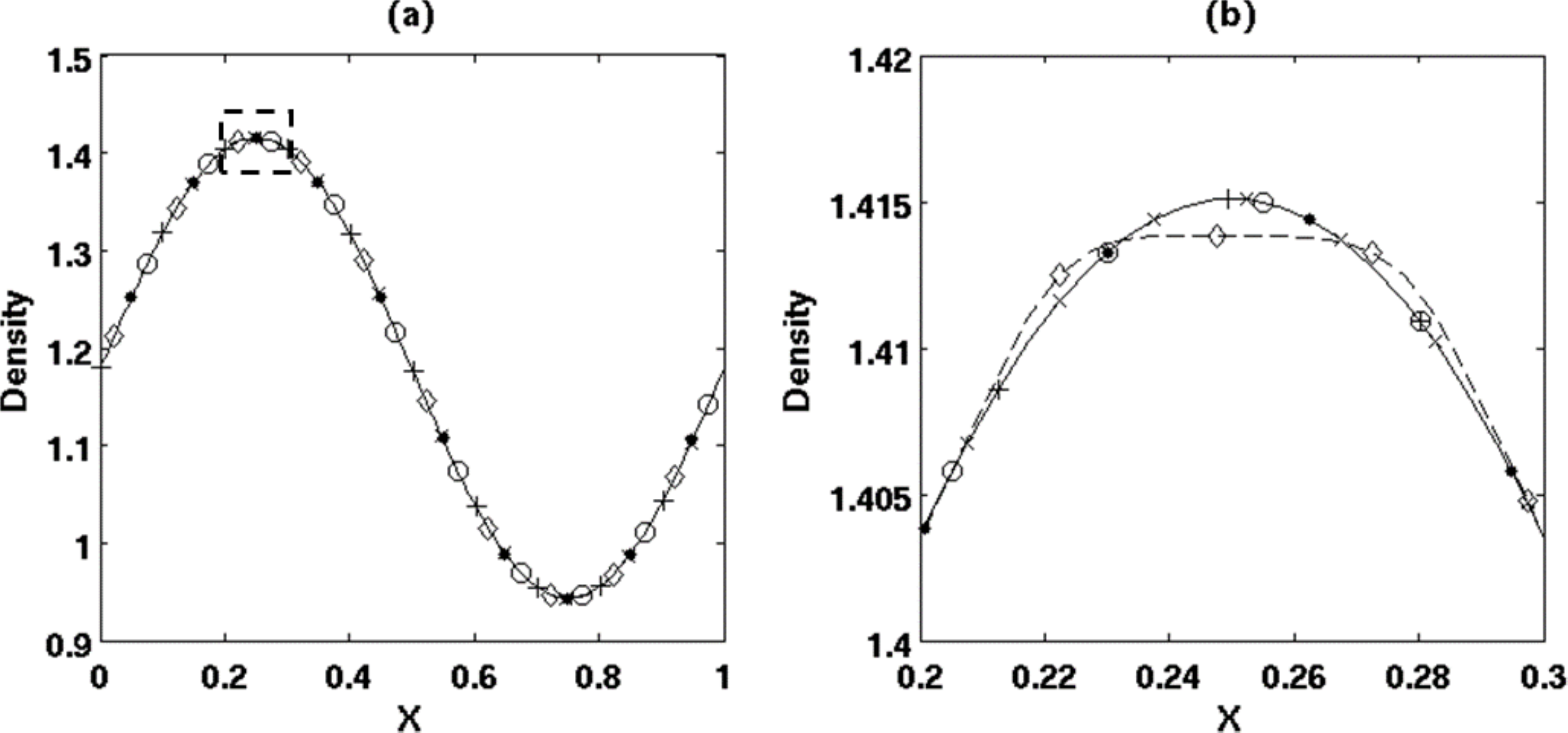} 
 \caption{Translating density sine wave. Density profile on $40$ SD elements in $x$ (equivalent to 200 FV cells for WENO\_AUSM and WENO\_AUSM\_TVD).  SD\_RUS, $\circ$;  SD\_AUSM, $\bullet$; WENO\_AUSM, $\times$; WENO\_AUSM\_TVD, $\Diamond$ with dashed line; HYBRID\_STATIC, +, and the exact solution; solid line.  }
 \label{sine_density}
\end{center}
\end{figure}

\section{Order of Accuracy}
\label{order_accuracy}
In this section, the order of accuracy of the $5^{th}$-order SD method, the $5^{th}$-order WENO scheme and the hybrid method with the $5^{th}$-order SD method and the $5^{th}$-order WENO scheme is examined.  For the SD method, two Riemann solvers are tested, namely the Rusanov solver  and the AUSM$^+-$up scheme. The Rusanov solver is commonly used in the SD method. The AUSM$^+-$up scheme is chosen for its simple yet robust approach covering very low to high speed flows, and it is known to be less prone to the carbuncle effect.  Each of these combinations is denoted by SD\_RUS and SD\_AUSM, respectively.  Based on the accuracy study with SD\_RUS and SD\_AUSM,  the AUSM$^+-$up scheme is chosen for the Riemann solver in the WENO scheme, denoted by WENO\_AUSM. Since the WENO scheme is primarily used as the shock-capturing scheme, a limiter is needed in some cases with strong shock. Therefore, the TVD limiter, Eq. ({\ref{tvd_lim}), is employed, denoted by WENO\_AUSM\_TVD, and its effect on the WENO solution is examined.  For the hybrid method, the AUSM$^+-$up scheme is used in both SD method and WENO scheme, and it is denoted by HYBRID\_STATIC as the grid is setup in static combination of SD elements and FV elements as shown in Fig. \ref{sine_vort_grid}. In all cases, the 3$^{rd}$-order strong stability preserving Runge-Kutta scheme is used. For the study of the order of accuracy in space, the time step should be sufficiently small in order for the spatial error to be dominant over the temporal error.  The time step size, thus the CFL number, is estimated based on the grid size of the  FV subcells in the FV element and  $\Delta t \approx \Delta x ^{n_x/n_t}$, where $n_x$ is the spatial order of accuracy and $n_t$ the temporal order of accuracy.  Based on this estimation, the time step is set sufficiently small.  The standard tests of translating density sine wave and translating vortex are used for accuracy study.  One dimensional arrangement of hexahedral elements is used in these test cases.
 
\subsection{Translating density sine wave}
First, the error and the order of accuracy of the $5^{th}$-order SD method, the $5^{th}$-order WENO scheme and the hybrid method are examined in 1D by a translating density sine wave.  The physical domain in $x$ is  $[0, 1]$ m.  The physical domain is discretized with 20, 40 and 80 base hexahedral elements in $x$ direction. For SD\_AUSM and SD\_RUS, the 5$^{th}$-order SD element is used in all base hexahedral elements.  For WENO\_AUSM and WENO\_AUSM\_TVD, the FV element with $5^3$ uniform FV subcells is embedded in all base hexahedral elements, resulting in the same solution DOF as in the $5^{th}$-order SD element and the resolution of 100, 200 and 400 FV subcells in $x$ direction. For HYBRID\_STATIC case, the physical domain is discretized with 20, 40 and 80 base hexahedral elements, and the FV element with $5^3$ FV subcells is statically embedded in the base hexahedral elements for $x > 0.5$ and the $5^{th}$-order SD elements for $x < 0.5$.  The density  profile is initialized  with a sine function  by 
\begin{equation}
\rho = \rho_0 ( 1 + 0.2 \sin (2 \pi x)),
\end{equation}
where $\rho_0$ is the reference density, $\rho_0 = 1.179$ kg/m$^3$.  The pressure is set constant at 1 atm, the temperature is computed from the  equation of state, $T = p/ \rho R$, with $R=288.18$ J/kg K. The density sine wave is translated with the velocity, $u=100$ m/s. The time step is set to $\Delta t = 5 \times 10^{-7}$ in all cases, such that the spatial error is dominant. Based on the grid size of FV subcells, the CFL numbers are 0.005, 0.01 and 0.02 for the spatial resolution of 100, 200 and 400 FV subcells, respectively.  Periodic boundary conditions are used in all directions.  The initial density profile in the mesh configuration for HYBRID\_STATIC case is depicted in Fig. \ref{sine_vort_grid}(a).  
 
The density sine wave is translated for two periods of a cycle.  Then, the density profile is compared with the exact solution.  The error norm and the order of accuracy of various cases are listed in Table \ref{sine_error}. It shows that the designed order of accuracy is nearly reached with SD\_AUSM, SD\_RUS and WENO\_AUSM. It should be noted that the error norm is smallest in SD\_AUSM among these three cases.  In the case of WENO\_AUSM\_TVD, the order of accuracy is much reduced to around $2^{nd}$ order. This is due to the cutoff behavior of statically enforced TVD limiter clipping local extrema, degrading the scheme to lower order. The HYBRID\_STATIC case also shows the degradation of the order to around $3^{rd}$ order.  It should be noted that there are two additional sources of the numerical errors in the hybrid method besides the truncation error associated with the numerical schemes employed. One is the interpolation of solutions between the set of Gauss points and  the set of FV subcell solution points in the SD element neighboring the FV element, as depicted in Fig. \ref{weno_interp}, and the other is the mortar projection process, interpolating the conservative variable from the  FV subcells to the mortar interface and  interpolating the computed common flux  back to the  FV subcells, as shown in Fig. \ref{coupling}.  These interpolation errors can be reduced by using more FV subcells in the FV element.  As suggested in Dumbser {\it et al.} \cite{dumbser2014}, the optimal number of  FV subcells would be $2N+1$ per spatial direction, where $N$ is the order of solution polynomial in the SD element.  Thus, using $5^3$ FV subcells to keep the same solution DOF as in the  $5^{th}$-order SD element is suboptimal and introduces relatively larger truncation error in the interpolation steps.  However, it should be noted that in this 1D test case, the mortar projection step does not play role, thus the interpolation error is mainly resulted from the interpolation of solutions between two sets of different solution points inside the SD element, leading to the reduction of the order. Although the order of accuracy in HYBRID\_STATIC is reduced, it should be pointed out that the magnitude of the error norm by HYBRID\_STATIC is much smaller than the error norm by WENO\_AUSM\_TVD and only an order of magnitude larger than the error norm by WENO\_AUSM as the mesh is refined. The density profiles of various cases are shown in Fig. \ref{sine_density}.  It shows that the density profiles of various cases are in good agreement with the exact solution except the WENO\_AUSM\_TVD case.  The WENO\_AUSM\_TVD case is depicted by dashed line with diamond mark in Fig. \ref{sine_density}(b), and it shows the cutoff behavior of the TVD limiter at local extrema.  

\begin{table}[t]
\begin{center}
\caption{Translating vortex: error and order of accuracy.}
\begin{tabular} { c|c|c c|c c|c c }
\hline \hline
               & Mesh                  &$L_1$ error & Order& $L_2$ error& Order& $L_{\infty}$ error &  Order \\ \hline 
\multirow{3}{*}
{SD\_AUSM}     & 20$\times$20$\times$1 & 2.5311e-06 & -    & 5.1440e-06 & -    & 6.2432e-05 & -  \\ 
               & 40$\times$40$\times$1 & 9.1350e-08 & 4.79 & 1.6210e-07 & 4.99 & 2.4500e-06 & 4.67 \\   
               & 80$\times$80$\times$1 & 2.5259e-09 & 5.18 & 5.3222e-09 & 4.93 & 8.9683e-08 & 4.77 \\ \hline 
\multirow{3}{*}
{WENO\_AUSM}   & 100$\times$100$\times$1 & 1.1542e-04 & -    & 3.1883e-04 & -    & 2.5294e-03 & - \\   
 			   & 200$\times$200$\times$1 & 2.7858e-05 & 2.05 & 7.7676e-05 & 2.04 & 5.4659e-04 & 2.21 \\
               & 400$\times$400$\times$1 & 7.9988e-06 & 1.80 & 1.9391e-05 & 2.00 & 1.3792e-04 & 1.99 \\  \hline 
\multirow{3}{*}
{HYBRID\_STATIC } 
               & 20$\times$20$\times$1 & 1.2735e-04 & -    & 2.8108e-04 & -    & 2.4218e-03 & - \\ 
               & 40$\times$40$\times$1 & 3.0947e-05 & 2.04 & 6.7799e-05 & 2.05 & 5.0281e-04 & 2.27 \\   
 			   & 80$\times$80$\times$1 & 7.7032e-06 & 2.01 & 1.7005e-05 & 2.00 & 1.2438e-04 & 2.02 \\ \hline \hline
\end{tabular}
\label{vortex_error}
\end{center}
\end{table}

\subsection{Translating vortex}
Based on the error and the order of accuracy study in 1D translating density wave, SD\_AUSM, WENO\_AUSM and HYBRID\_STATIC cases are chosen for the order of accuracy test in 2D translating vortex. The physical domain in $x$ and $y$ is  $[0,10] \times [0,10]$.  The physical domain is discretized by $20 \times 20$, $40 \times 40$ and $80 \times 80$ base hexahedral elements.  For SD\_AUSM, the 5$^{th}$-order SD element is used in all base hexahedral elements.  For WENO\_AUSM, the FV element with $5^3$ FV subcells is embedded in all base hexahedral elements, resulting in the resolution of $100 \times 100$, $200 \times 200$ and $400 \times 400$ FV subcells in $x$ and $y$ directions constructing the same solution DOF as in SD\_AUSM.   For HYBRID\_STATIC case, the physical domain is discretized in the same way as in SD\_AUSM except that the FV element with $5^3$ FV subcells is embedded in the base hexahedral elements for $x > 5$. The mean flow conditions are $(\rho, u, v, w, p) = (1, 1, 0, 0, 1)$. The isotropic vortex is superimposed on the mean flow with the following perturbations.  

\begin{eqnarray}
\Delta u &=& \frac{ \epsilon}{2 \pi} e^{( 1- r^2)/2 } (5-y), \\ 
\Delta v &=& \frac{ \epsilon}{2 \pi} e^{( 1- r^2)/2 } (x-5), \\  
\Delta T &=& - \frac{(\gamma - 1) \epsilon^2} {8 \gamma \pi^2} e^{(1-r^2)}, \\
\Delta S &=& 0,
\end{eqnarray} 
where $r^2 = (x-5)^2 + (y-5)^2$, the vortex strength, $\epsilon = 5.0$, and the ratio of specific heats, $\gamma=1.4$.  Initially, the center of isotropic vortex is located at $x=5$ and $y=5$. The time step is set to $\Delta t = 2 \times 10^{-3}$ in all cases, such that the spatial error is dominant. Based on the grid size of  FV cells, the CFL numbers are 0.02, 0.04 and 0.08 for the spatial resolution of $100 \times 100$, $200 \times 200$ and $400 \times 400$ FV subcells. Periodic boundary conditions are used in all directions.  The initial density profile in the mesh configuration for HYBRID case is depicted in Fig. \ref{sine_vort_grid}(b).  

The vortex is translated for two periods of a cycle.  Then, the density profile is compared with the exact solution.  The error norm and the order of accuracy of three cases are listed in Table \ref{vortex_error}.  It shows that the designed order of accuracy is achieved only in SD\_AUSM while the order of accuracy is much reduced in WENO\_AUSM and HYBRID\_STATIC.  Examining the location of maximum error norms of density in WENO\_AUSM and HYBRID\_STATIC, the maximum error norms are observed in the vortex core region, where the gradients of velocity components are steepest and the pressure gradient is rapidly changing.   Then, these gradients may have introduced additional numerical dissipation by  the pressure and the velocity dissipation terms in mass and pressure fluxes, respectively, in the AUSM$^+-$up scheme at low Mach number while these dissipation terms are automatically turned off in the supersonic flow. Readers are referred to Liou \cite{liou2006ausm} for the discussions regarding numerical dissipation terms in mass and pressure fluxes.  Especially, the dissipation term by pressure gradient in the mass flux in the AUSM$^+-$up scheme seems to contribute additional numerical dissipation in this case while this term is zero in the translating density sine wave with the constant pressure. Thus, the order of accuracy may have been reduced  in WENO\_AUSM due to these numerical dissipation in mass and pressure fluxes in AUSM$^+-$up scheme. In HYBRID\_STATIC, there still exist interpolation errors as discussed in translating density sine wave. However, by examining and comparing the location and magnitude of error norms with the error norms in WENO\_AUSM and HYBRID\_STATIC, the numerical dissipation in mass and pressure fluxes in the AUSM$^+-$up scheme seems to be more dominant over the interpolation errors, thus resulting in about the same order of accuracy as WENO\_AUSM.  It is discussed in Scandaliato and Liou \cite{scandaliato2010} that the AUSM$^+-$up scheme has not been much tested with the high-order reconstruction scheme, and  some suppressed errors may start to appear when AUSM$^+-$up scheme is coupled with high-order scheme, as opposed to coupling the AUSM$^+-$up scheme with the low-order scheme.  In addition, it is also discussed in Hosangadi {\it et al.} \cite{hosangadi2012} that substantial dissipation is added to the pressure flux in AUSM$^+-$up for unsteady low Mach number flows. The issue with the numerical dissipation in the AUSM$^+-$up scheme with the high-order reconstruction needs to be further investigated, and it will be reported in the future work.

While the WENO reconstruction with the AUSM$^+-$up scheme is used in the statically embedded FV elements in the domain where $x > 5$ in HYBRID\_STATIC, it should be noted that the WENO reconstruction with the AUSM$^+-$up scheme is primarily employed  in the proposed hybrid method as a shock-capturing scheme dynamically applied only in the vicinity of the region where discontinuities  are present, thus the order of accuracy of the WENO with the AUSM$^+-$up schemes is degraded only in these regions. Therefore, it is expected that the degradation of the order of accuracy of the WENO reconstruction with the AUSM$^+-$up scheme does not affect the overall performance of the proposed hybrid method. (See Fig. 6.)

\section{Test Cases}
\label{test_cases}
In this section, the standard 1D and 2D test cases are presented using the proposed hybrid method and the $5^{th}$-order WENO scheme, denoted by HYBRID and WENO respectively, for comparison.  In HYBRID, the FV element with structured $5^3$ FV subcells is embedded dynamically only in the base hexahedral elements where the discontinuities are present and the $5^{th}$-order SD element in the smooth flow region away from the discontinuities. For WENO, the FV element with structured $5^3$ FV subcells are embedded in all base hexahedral elements, resulting in 5 times fine grid in all directions to have the same solution DOF as in the $5^{th}$-order SD element.  For all test cases presented, the AUSM$^{+}$-up scheme \cite{liou2006ausm} is employed as the Riemann solver.  For the test cases with a strong shock, the TVD limiter is used.  In all test cases, hexahedral elements are arranged in the $x$ direction (1D) or in the $x-y$ plane (2D) with one hexahedral element in $z$ direction.

\begin{figure}[t]
\begin{center} 
 \includegraphics[width=6in]{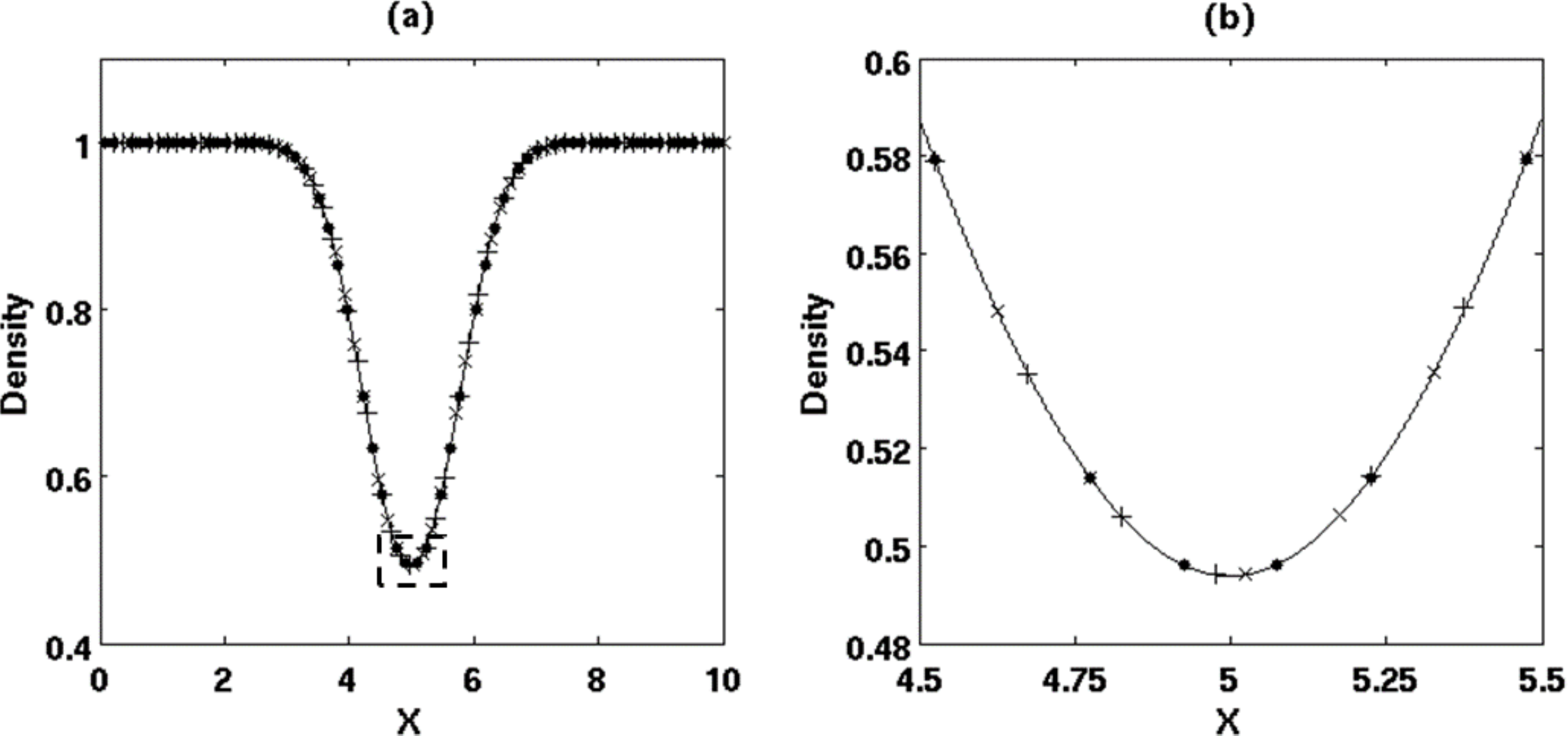}   
 \caption{Translating vortex. Density profile on $40 \times 40$ SD elements in $x$ and $y$ (equivalent to $200 \times 200$ FV subcells for WENO\_AUSM). SD\_AUSM, $\bullet$; WENO\_AUSM, $\times$; HYBRID\_STATIC, +, and the exact solution; solid line.}
 \label{tv_density}
 \end{center}
\end{figure}

\subsection{SOD shock tube problem}  
The SOD shock tube problem \cite{sod1981} is considered to examine the performance of shock-capturing capability by dynamically embedding the FV element in the base hexahedral elements following the shock and the contact surface. The computed solutions between HYBRID and WENO are compared. The physical domain in $x$ is $[0, 1]$ m and is discretized by 20, 40, 80 and 160 base hexahedral elements. The base grid results in 100, 200, 400 and 800 FV cells in $x$ direction for WENO. The initial condition is set  by
\begin{eqnarray}
(\rho, u, v, w, p) &=& (10 \rho_0, 0, 0, 0 , 10 p_0) \quad \mbox{for} \quad x < 0.5, \\
(\rho, u, v, w, p) &=& (\rho_0, 0, 0, 0 , p_0) \quad \quad \quad \mbox{for} \quad x \geq 0.5,
\end{eqnarray} 
where $\rho_0 = 1.179$ kg/m$^3$ and  $p_0 = 1$ atm. The specific gas constant, $R$, is 288.18 J/kg K.
The density and pressure discontinuities in the SD elements are initially detected, and the FV element with $5^3$ FV subcells is embedded in place of the SD element at the locations of discontinuities. Then, the initial conditions in the SD element with discontinuities are interpolated onto the FV subcell solution points.  As the solution is advanced in time with constant $\Delta t = 1 \times 10^{-6}$, the contact surface and the shock are dynamically detected. The discontinuity detector is used on both density and pressure in order to capture the contact surface and the shock with $\epsilon_s = 0.01$.

\begin{figure}[t!]
\begin{center}
 \includegraphics[width=6in]{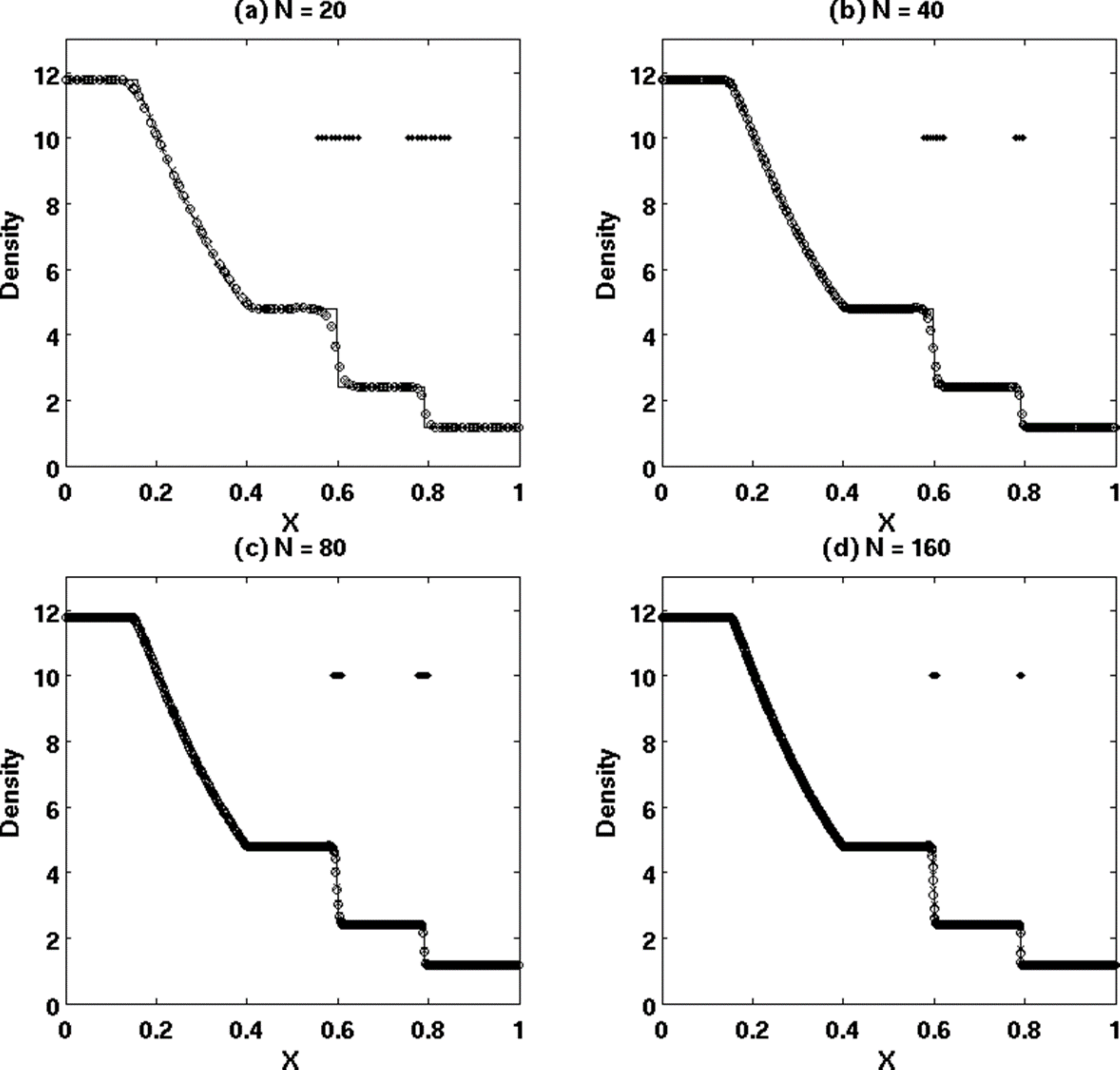} 
 \caption{SOD shock tube problem. Density profiles on different grid sizes with $N$ denoting the number of base hexahedral elements. HYBRID, $\circ$; WENO, $\times$ and the exact solution, solid line.  Horizontal dots indicate the locations of FV subcell solution points.}
 \label{shock_tube_1}
\end{center}
\end{figure}

\begin{figure}[t!]
\begin{center}
 \includegraphics[width=6in]{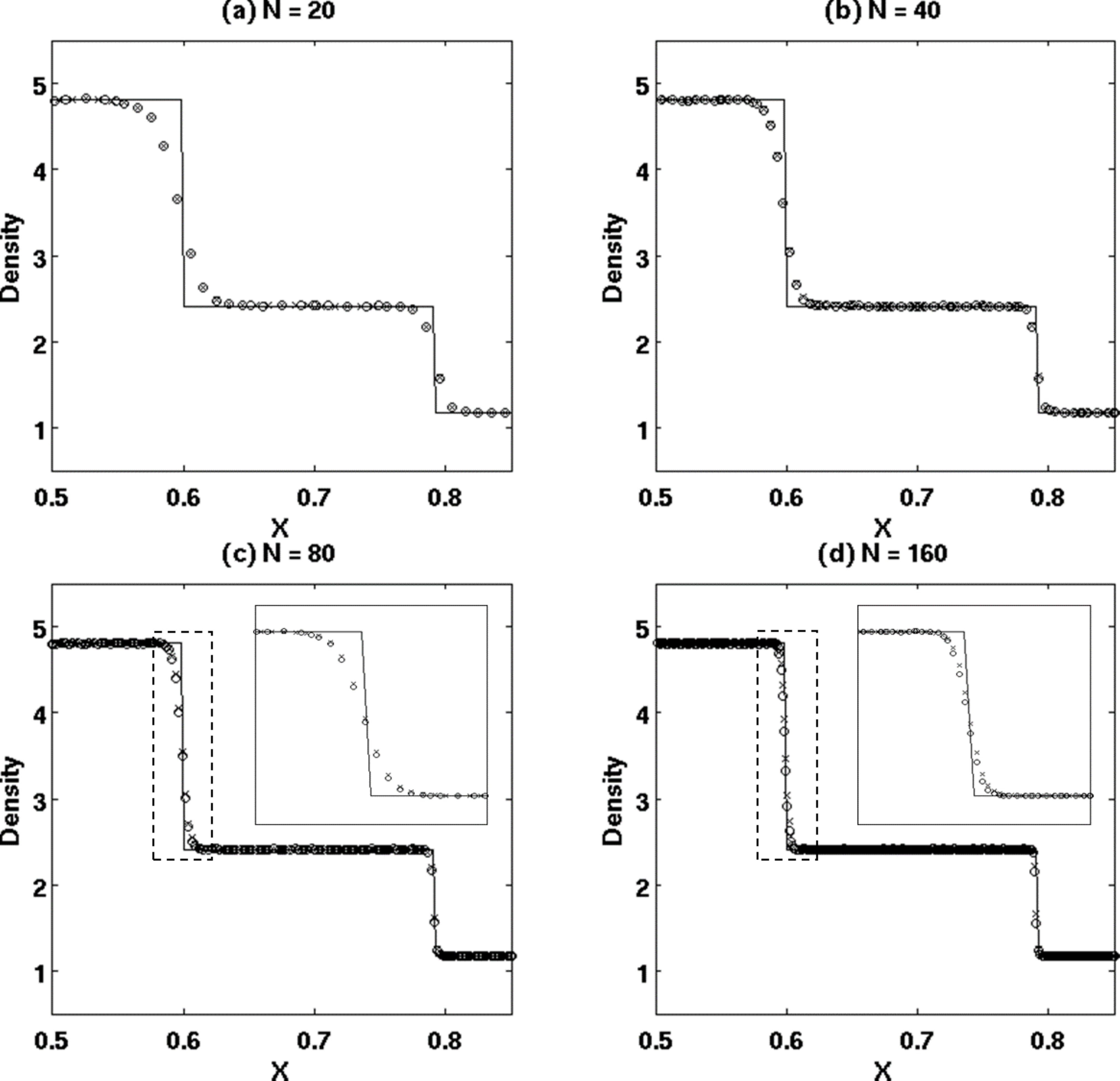} 
 \caption{SOD shock tube problem. The close-up view around the contact surface and the shock.    HYBRID, $\circ$; WENO, $\times$ and the exact solution, solid line.}
 \label{shock_tube_2} 
 \end{center}
\end{figure}

The computed density profiles are shown in Fig. \ref{shock_tube_1} for various grid sizes with $N$ denoting the number of the base hexahedral elements. Figure \ref{shock_tube_2} shows the close-up views near the contact surface and the shock. The exact solution is also shown in these figures for comparison. The horizontal dots in Fig. \ref{shock_tube_1} indicate the locations of the FV subcell solution points.  It can be seen that the region where the FV element is embedded in HYBRID becomes smaller close to the discontinuities and that the discontinuities are captured at the FV subcell resolution within one or two embedded FV elements as the mesh is refined.  The density profiles computed by both HYBRID and WENO compare well with the exact solution for the given spatial resolutions. The density profiles show that HYBRID with coarser mesh produces comparable results as WENO with 5 times finer grid for the same solution DOF. For example, the hybrid method with 160 hexahedral elements produces very similar results as the $5^{th}$-order WENO scheme on 800 FV cells. The similar conclusion is also discussed in Dumbser {\it et al.}  \cite{dumbser2014} obtaining an accurate solution with a coarse mesh.  Overall, by using the current hybrid method, the contact surface and shocks are captured at the sub-element resolution compared to other approaches, such as using artificial viscosity or limiter in the discontinuous high-order methods, where the discontinuities are usually captured in a few elements. Comparing the density profile on the 20 base hexahedral elements, Figs. \ref{shock_tube_1}(a) and \ref{shock_tube_2}(a), with the results by Dumbser {\it et al.}  \cite{dumbser2014} for the same number of base element, it is noted that the shock and the contact discontinuity show slightly more diffusive profile than the one shown in Dumbser {\it et al.}.  It is because the number of subcells used in Dumbser {\it et al.} is $2N+1 = 19$ with $N=9$ while the number of subcells in this study is 5.  Therefore, the discontinuity in Dumbser {\it et al.} is much sharply resolved.  It is expected that when the same number of subcells are used as in Dumbser {\it et al.}, the results would be almost the same.

\begin{figure}[t!]
\begin{center}
 \includegraphics[width=6in]{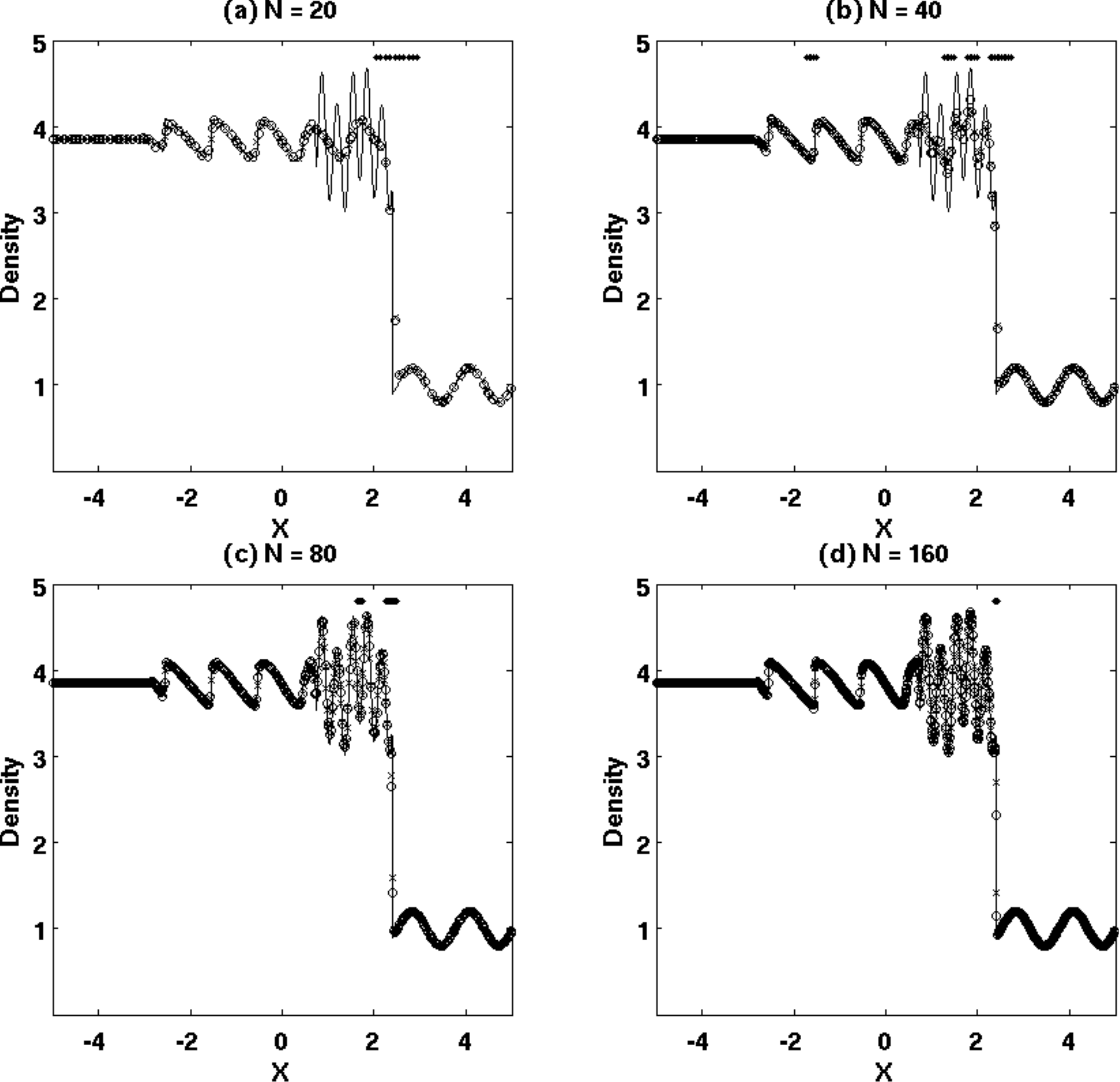} 
 \caption{Shu-Osher problem. Density profiles on different grid sizes with $N$ denoting  the number of base hexahedral elements.   HYBRID, $\circ$; WENO, $\times$ and the numerical "exact" solution, solid line.  Horizontal dots indicate the locations of FV subcell solution points}
 \label{entropy_1}
\end{center}
\end{figure}

\begin{figure}[t!]
\begin{center}
 \includegraphics[width=6in]{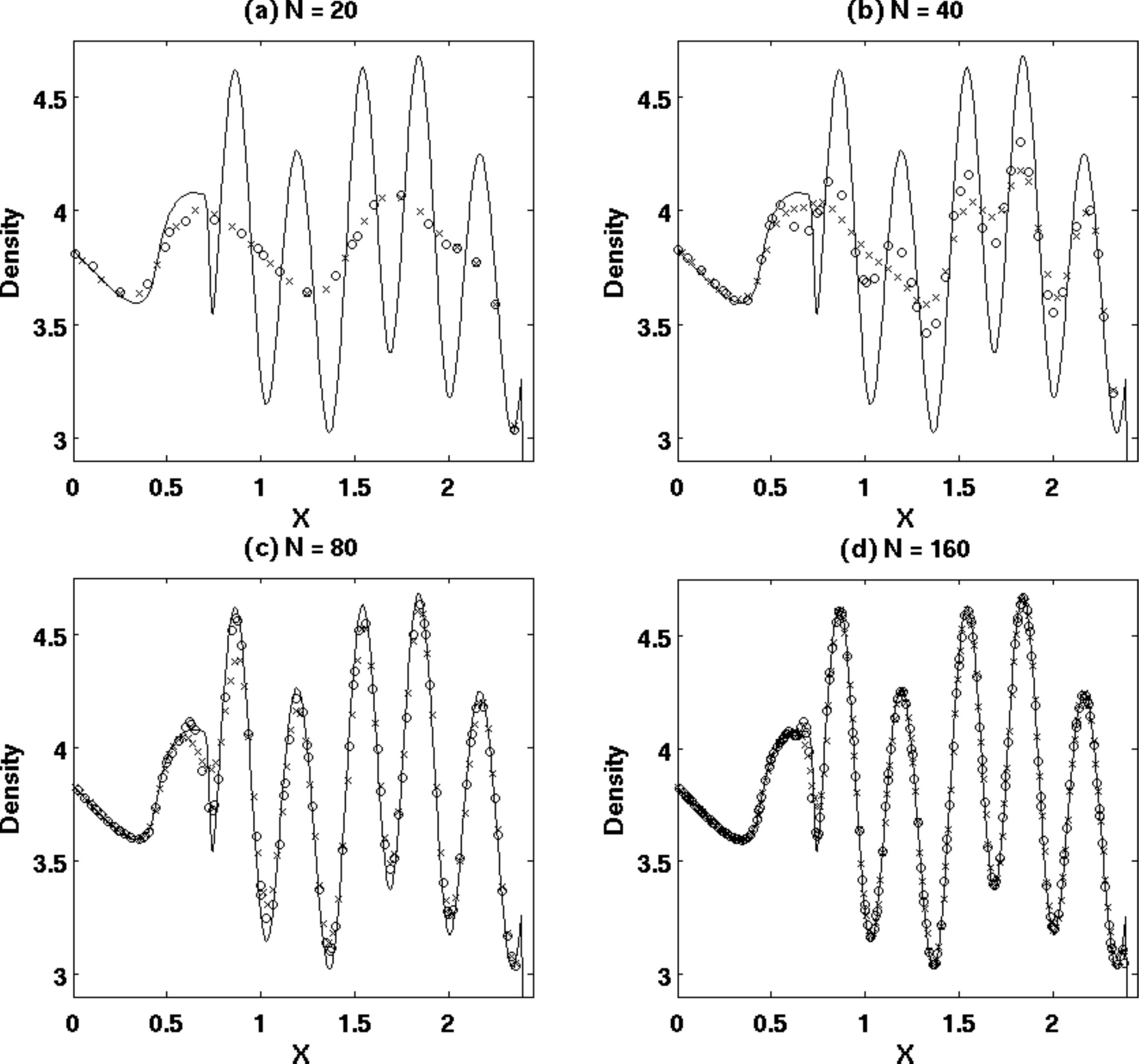} 
 \caption{Shu-Osher problem. The close-up view around the compress entropy waves behind the shock.   HYBRID, $\circ$; WENO, $\times$ and the numerical "exact" solution, solid line.}
 \label{entropy_2}
 \end{center}
\end{figure}
    
\subsection{Shu-Osher problem}  
In this example, the interaction between propagating shock and perturbed entropy waves is simulated \cite{shu1989}.  The physical domain in $x$ is $[-5, 5]$. The domain is discretized with 20, 40, 80 and 160 base hexahedral elements, which results in  100, 200, 400 and 800 FV cells in $x$ direction for WENO.  The nondimensional initial condition is given by
\begin{eqnarray}
(\rho, u, v, w, p) &=& (3.857143, 2.269369, 0, 0, 10.33333) \quad \mbox{for} \quad x < -4, \\
(\rho, u, v, w, p) &=& (1+0.2\sin(5x), 0, 0, 0, 1) \quad \quad \quad \quad \quad
~ \mbox{for} \quad x \geq -4.
\end{eqnarray} 
with the ratio of specific heats, $\gamma=1.4$.

The initial shock is detected in SD elements, and the FV element is embedded in these elements. The solution is advanced with constant  $\Delta t=1 \times 10^{-4}$. The FV element is dynamically embedded replacing the SD element in the hexahedral elements where the steep density gradient or the shock is detected with $\epsilon_s = 0.02$.  The numerical "exact" solution is obtained by the $5^{th}$-order WENO scheme employed in this study on a 3000 FV cells.

\begin{figure}[t!]
\begin{center}
 \includegraphics[width=5in]{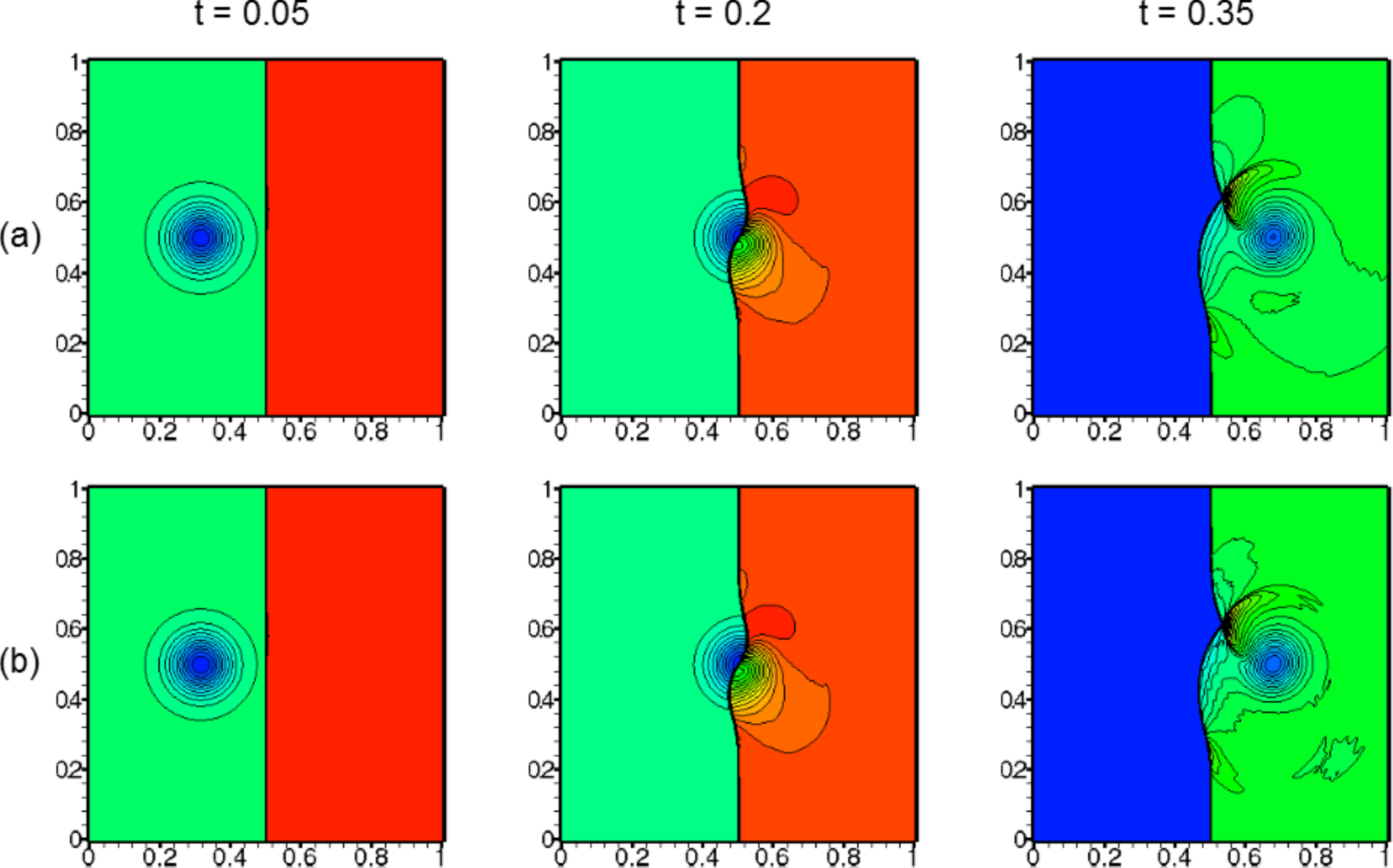} 
 \caption{Weak shock/vortex interaction. Density contours with  20 contour lines from minimum to maximum using (a)  WENO  and (b) HYBRID at $t=0.05$, $t=0.2$ and $t=0.35$.}
 \label{shock_vortex_1}
 \end{center}
\end{figure}

\begin{figure}[t!]
\begin{center}
 \includegraphics[width=5in]{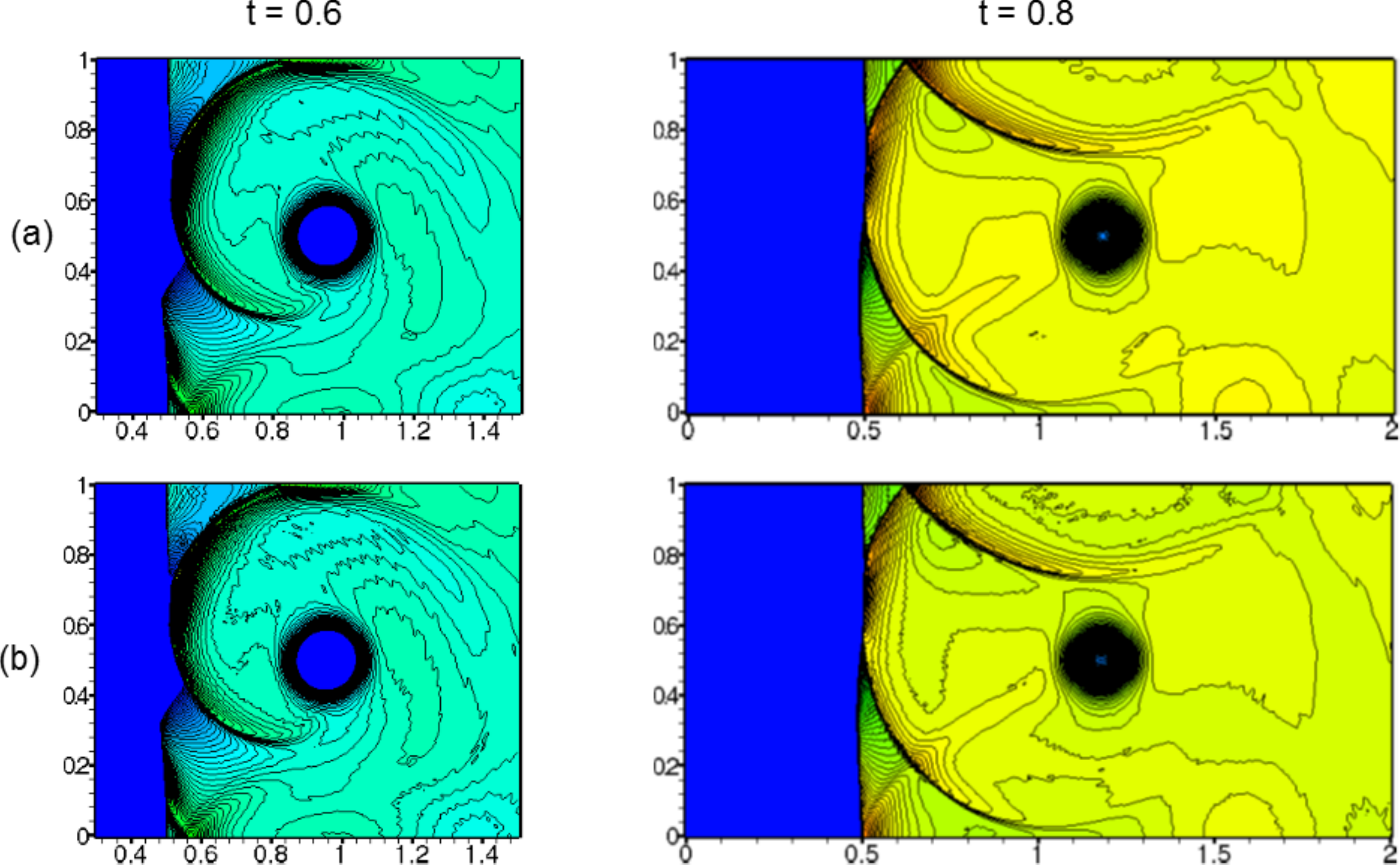} 
 \caption{Weak shock/vortex interaction. Density contours  using (a) WENO and (b) HYBRID at $t=0.6$ and $t=0.8$. At $t=0.6$, 90 contour lines from 1.19 to 1.37 and at $t=0.8$, 90 contour lines from minimum to maximum.}
 \label{shock_vortex_2}
 \end{center}
\end{figure}

The simulation is run until $t=1.8$. The computed density profiles as well as the numerical "exact" solution are shown in Fig. \ref{entropy_1} for various grid sizes with $N$ denoting the number of base hexahedral elements. The close-up views are shown in Fig. \ref{entropy_2}. In Fig. \ref{entropy_1}, the horizontal dots show the locations of FV subcell solution points. The computed results show that the compressed entropy waves behind the shock are not well capture by both HYBRID and WENO at low resolutions, e.g. $N=20$ and $N=40$. But, it can be seen that  HYBRID gives a slightly better result with $N=40$ than WENO.  When the grid is refined further to $N=80$ and $N=160$, both HYBRID and WENO, having the same solution DOF, are able to capture the compressed entropy waves better showing good agreement with the numerical "exact" solution. Comparing with the results by Dumbser {\it et al.} \cite{dumbser2014}, the similar conclusion can be drawn as discussed in the SOD shock tube problem.

\begin{figure}[t!]
\begin{center}
 \includegraphics[width=5in]{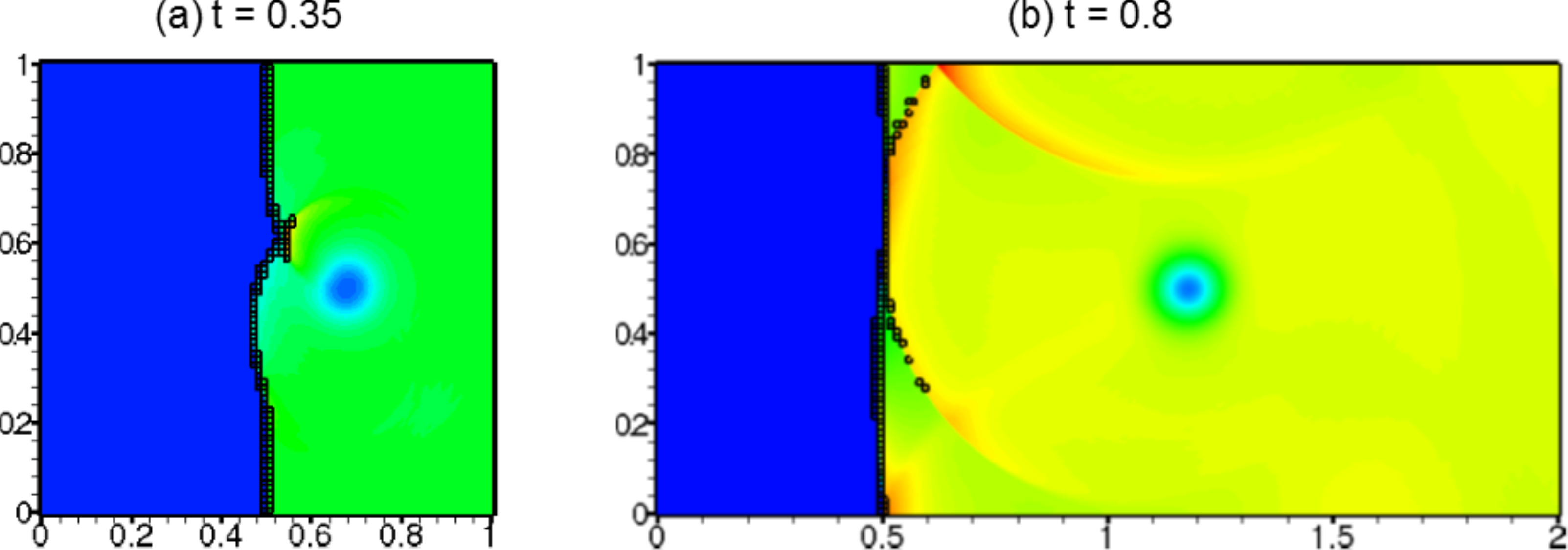} 
 \caption{Weak shock/vortex interaction. The black squares show the locations of embedded FV elements at (a) $t=0.35$ and (b) $t=0.8$. }
 \label{shock_vortex_3}
 \end{center}
\end{figure}

\subsection{Weak shock/vortex interaction}
This test simulates the interaction of a moving vortex with a stationary weak shock \cite{jiang1996}. The physical domain in $x$ and $y$ is $[0,2] \times [0,1]$, and the domain is discretized by uniform  $80 \times 40$ base hexahedral elements, which gives $400 \times 200$ FV cells for WENO. At $x=0.5$, a stationary shock with Mach 1.1 positioned, where the left state condition is $(\rho, u, v, w, p) = (1, 1.1 \sqrt{\gamma}, 0, 0, 1)$ with the ratio of specific heats, $\gamma=1.4$. A small vortex is superposed to the flow left of the shock at $(x_c, y_c) = (0.25, 0.5)$. The vortex is defined by perturbations to the mean flow defined by the left sate of the shock, and the perturbations are given by 
\begin{eqnarray}
\Delta u &=& \epsilon \tau e^{ \alpha ( 1- \tau^2) } \sin \theta, \\ 
\Delta v &=& -\epsilon \tau e^{ \alpha ( 1- \tau^2) } \cos \theta, \\
\Delta T &=& - \frac{(\gamma - 1) \epsilon^2 e^{2 \alpha (1-\tau^2)}  }{4 \alpha \gamma}, \\
\Delta S &=& 0,
\end{eqnarray} 
where $\tau = r/r_c$ and $r=\sqrt{(x-x_c)^2+(y-y_c)^2}$. $\epsilon$ is the strength of the vortex, $\alpha$ is the decay rate of vortex. $r_c$ is the critical radius for which the vortex has the maximum strength. In this study, $\epsilon=0.3$, $r_c=0.05$ and $\alpha=0.24$ are used. The reflective boundary conditions are used on the upper and lower boundary. The simulation is carried out with constant $\Delta t=1 \times 10^{-4}$. The discontinuity detector is used only on the pressure with $\epsilon_s = 0.01$.

The density contours at $t=0.05$, $t=0.2$ and $t=0.35$ are shown in Fig. \ref{shock_vortex_1}(a) for WENO and Fig. \ref{shock_vortex_1}(b) for HYBRID, and Fig. \ref{shock_vortex_2} (a) and (b) shows the density contours at $t=0.6$ and $t=0.8$ for WENO and HYBRID, respectively.   It can be seen that the shock is sharply captured in both  WENO and HYBRID by the $5^{th}$-order WENO scheme at the same resolution. The moving vortex is captured by the SD method in HYBRID.  It is shown that the shape of vortex after interacting with shock is still well defined in the simulations by both methods. At $t=0.8$, one of the shock bifurcations is reached the top boundary and reflected.  Overall, the simulations by both methods show very similar results with sharp shock-capturing and vortex preservation after the interaction with the shock.   Figure \ref{shock_vortex_3} shows the base hexahedral elements where the FV elements are embedded at $t=0.35$ and $t=0.8$, indicated by the black squares. It can be seen that the shock is captured at the FV subcell resolution within one or two hexahedral elements, of which it depends on the location of shock within the hexahedral element.  The weak reflected shock in HYBRID is captured by the SD method without severe oscillations as shown in Fig. \ref{shock_vortex_3}(b).

\subsection{Double Mach reflection}
The double Mach reflection is commonly employed to test high resolution schemes \cite{woodward1984}.  The physical domain in $x$ and $y$ is $[0,4] \times [0,1]$. This test case is carried out on two mesh sizes, $80 \times 20$ and $160 \times 40$ base hexahedral elements, which gives $400 \times 100$ and $800 \times 200$ FV cells for WENO. A right-moving Mach 10 shock is positioned at $x_0=1/6$ at an angle of $\pi/3$ with $x$ axis. The post-shock condition is ($\rho$, $u$, $v$, $w$, $p$) = (8, 8.249, -0.07539, 0, 116.5) and the ratio of specific heats, $\gamma=1.4$. In the region from $x=0$ to $x=1/6$ at the bottom boundary, the exact post-shock condition is imposed.  A reflective boundary condition is used for the remaining part at the bottom boundary.  On the top boundary, the exact motion of initial Mach 10 shock is imposed with the position of the shock wave at time $t$ at $y=1$  given by
\begin{equation}
s(t) = x_0 + \frac{1+20t}{\sqrt{3}}.
\end{equation}

The left boundary condition is set by the initial post shock condition, and the right boundary condition is set by the zero gradient condition. The simulation is run with constant $\Delta t=1 \times 10^{-5}$ until $t=0.2$. The TVD limiter is used, and the discontinuity detector with $\epsilon_s = 0.02$ is used on the pressure only.

The density contours by WENO and HYBRID on two mesh sizes are shown in Fig. \ref{dmr_1}(a) and (b), respectively.  The close-up views near Mach stems and jet formed are shown in Fig. \ref{dmr_2}.  At low resolution, both WENO and HYBRID produce comparable shock and jet structure. On the other hand, the jet structure is better resolved with HYBRID on the grid by $160 \times 40$ base hexahedral elements compared to the one by WENO on $800 \times 200$ FV cells for the same solution DOF.  It should be mentioned that the results obtained by the currently implemented  $5^{th}$-order WENO scheme on the $800 \times 200$ grid is comparable to the results shown in \cite{xu2005} and \cite{shi2003}, but the results by the $5^{th}$-order WENO scheme can certainly be improved by using more advanced variants of the WENO scheme, e.g., WENO-Z \citep{borges2008}.  In Fig. \ref{dmr_1}(c),  the location of base hexahedral elements where the FV element is embedded is shown, and it can be seen that the FV elements are embedded following shocks and Mach stems.  From this test case, it is noted that HYBRID performs better capturing dynamically developing smooth but complex flow structures than WENO  due to much less numerical dissipation in SD method compared with the WENO scheme.  In addition, it should be noted that  numerical dissipation from the TVD limiter used in the embedded FV elements (e.g. shock, Mach stem regions) in HYBRID does not propagate much into the regions with SD elements (e.g. jet region) away from discontinuities, such that the numerical dissipation from the TVD limiter does not pollute severely the small scale flow structure resolved by the low-dissipative SD method.   Comparing with the results in Sonntag and Munz \cite{sonntag2014} and Dumbser {\it et al.} \cite{dumbser2014}, where the base mesh resolution is 5 times and about 2 times finer, respectively, than the current fine mesh resolution, $160 \times 40$ base hexahedral elements, it is seen that the jet structures are well resolved, and  Kelvin–Helmholtz instabilities along the slip line are visible at the current resolution employed in the hybrid method.

\begin{figure}[t!]
\begin{center}
 \includegraphics[width=5.9in]{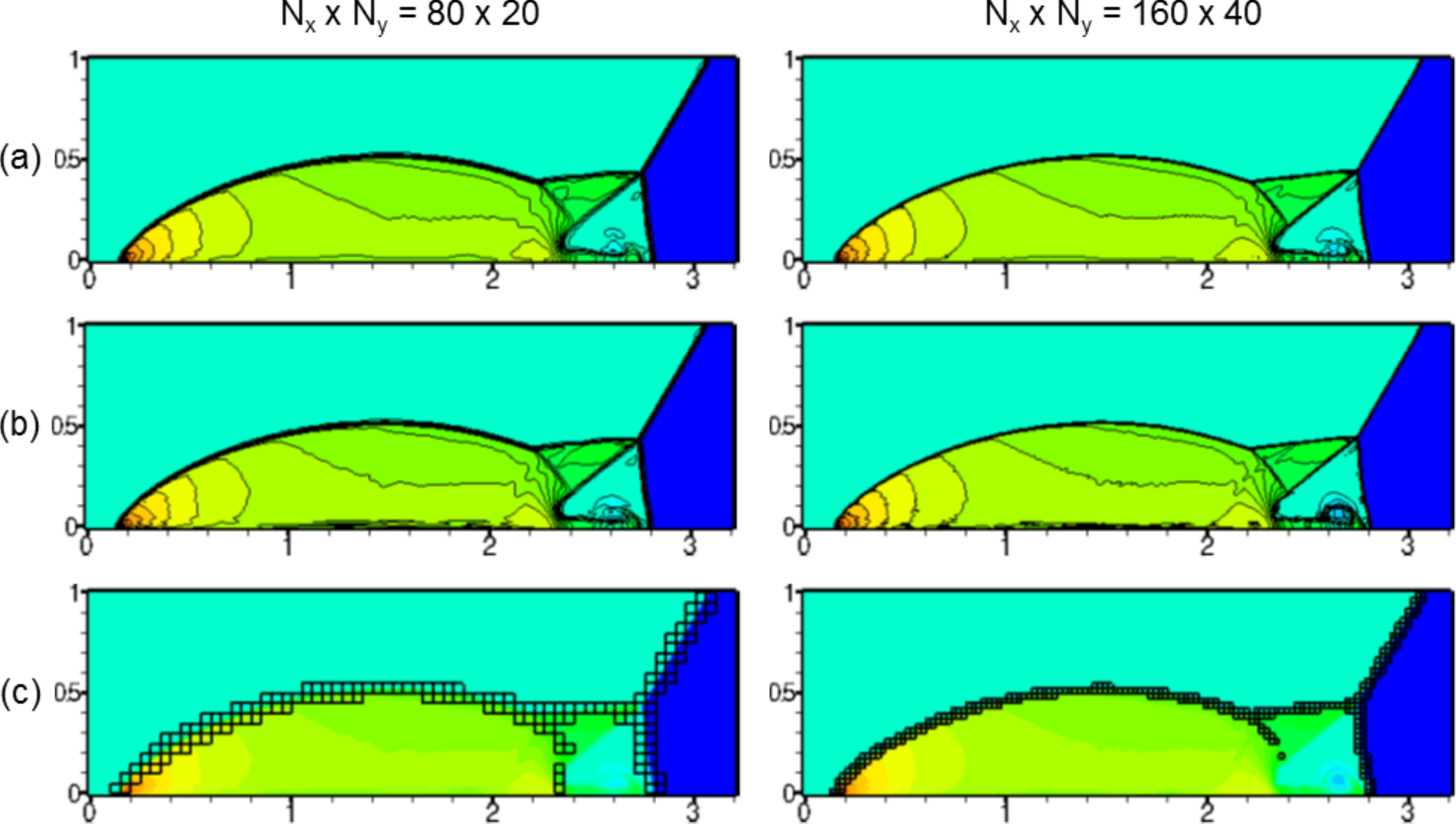} 
 \caption{Double Mach reflection. Density contours with 30 contour lines using (a) WENO and (b) HYBRID on two different grids with $N_x$ and $N_y$ denoting the number of base hexahedral elements in $x$ and $y$. (c) The locations of embedded FV elements are indicated by the black squares. }
 \label{dmr_1}
 \end{center}
\end{figure}

\begin{figure}[t!]
\begin{center}
 \includegraphics[width=4.9in]{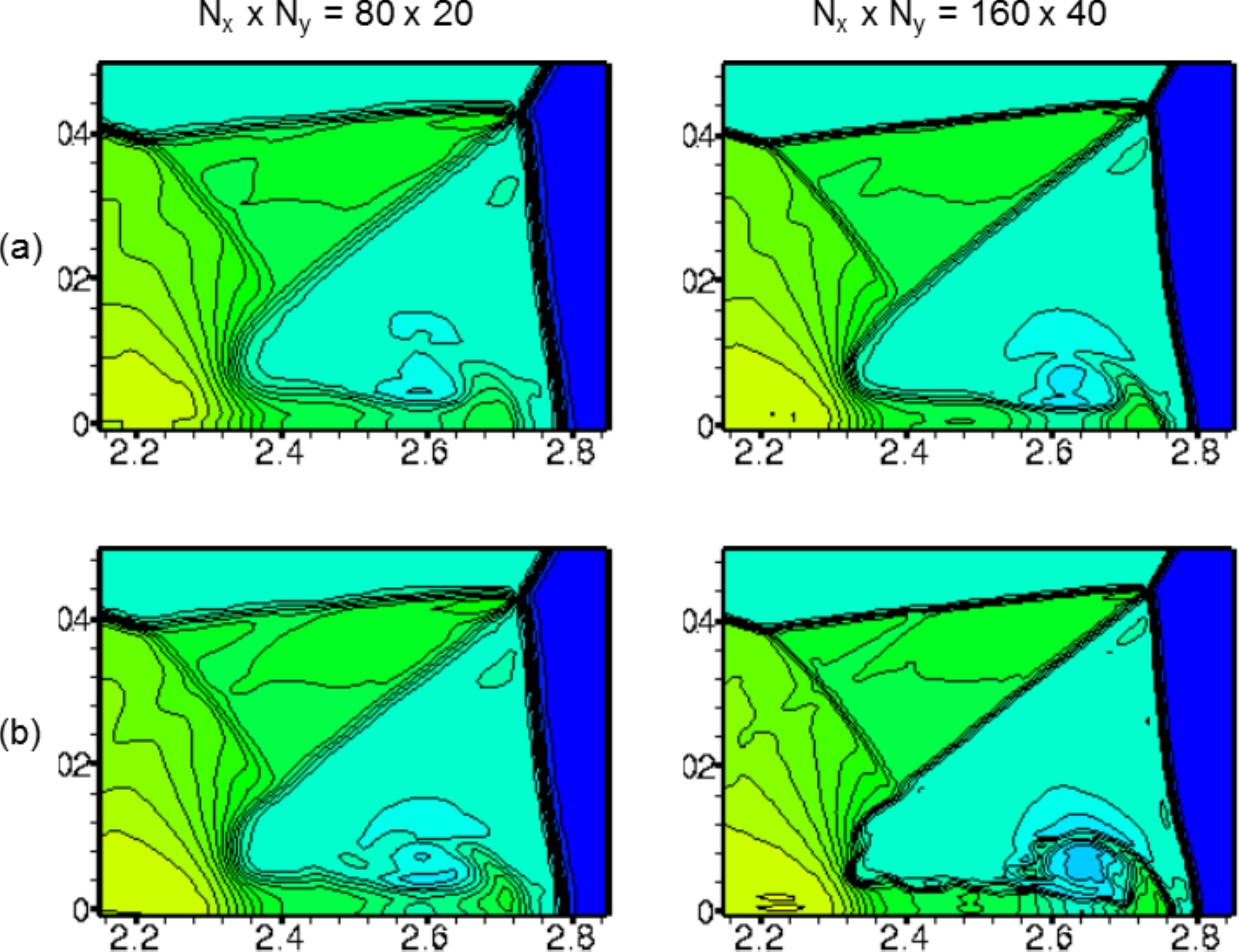} 
 \caption{Double Mach reflection. The close-up view of density contours near the double Mach stems using (a) WENO and (b) HYBRID.}
 \label{dmr_2}
 \end{center}
\end{figure}

\subsection{Inviscid supersonic flow past cylinder}
The inviscid supersonic flow past cylinder case is employed for testing the hybrid method on a curvilinear grid. The elliptic grid is constructed with $a=3$ m and $b=6$ m, where $a$ and $b$ are the minor and major axis of elliptic equation. The domain is discretized with $40 \times 30$ base hexahedral elements, resulting in  $200 \times 150$ FV cells for WENO.  The supersonic flow of air at Mach number 2 passes over the cylinder with the radius of 1 m.  The flow fields are initialized with $\rho=1.179$ kg/m$^3$, $p=1$ atm, and $T=298.15$ K and $u=693.31$ m/s.  The specific gas constant, $R$, is set to 288.18 J/kg K. On the cylinder wall, the inviscid wall condition is imposed.   The simulation is run with $\mbox{CFL} = 0.3$ until the solution is converged. The TVD limiter is employed, and the discontinuity detector is used on the pressure with $\epsilon_s =0.05$.
 
The pressure and density contours are shown in Fig. \ref{shock_cylinder_1}(a) and (b), respectively.  It can be seen that the shock is essentially captured at the resolution of  FV subcells in HYBRID, which is the same resolution in WENO. Thus, the simulation results show comparable solution between HYBRID and WENO.  Figure \ref{shock_cylinder_grid}(a) shows the location of base hexahedral elements where the FV element is embedded.  It should be noted that the FV element is statically embedded on the cylinder wall in order to approximate cylinder wall at the resolution of FV subcells, as shown in Fig. \ref{shock_cylinder_grid}(b). The issue with the curved boundary representation for the high-order scheme has been addressed in \cite{gao2010}.

\begin{figure}[t!]
\begin{center}
 \includegraphics[width=5.7in]{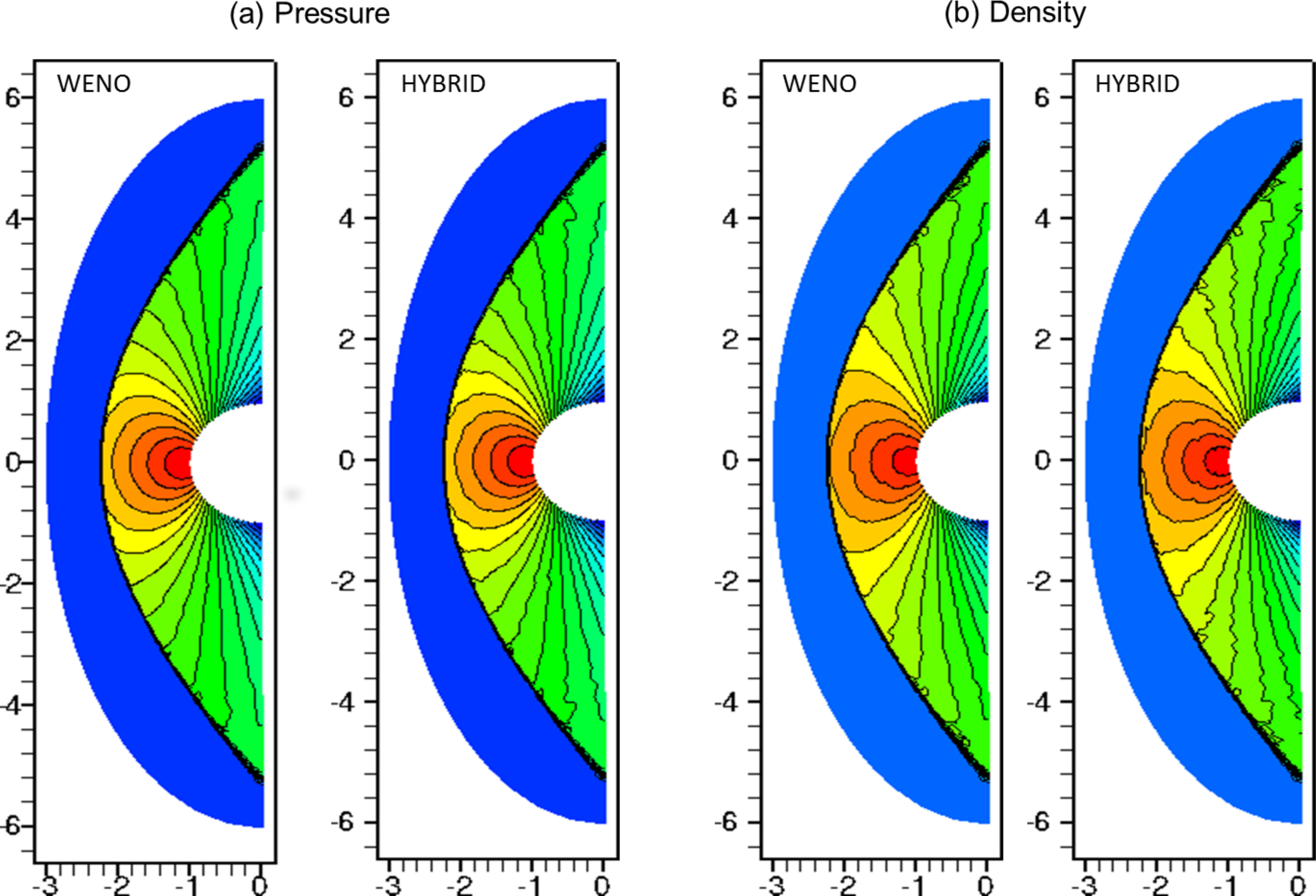} 
 \caption{Mach 2 flow over cylinder. Contours with 20 contour lines of (a) pressure and (b) density using  WENO and HYBRID.}
 \label{shock_cylinder_1}
 \end{center}
\end{figure}

\begin{figure}
\end{figure}

\begin{figure}[t!]
\begin{center}
 \includegraphics[width=5.7in]{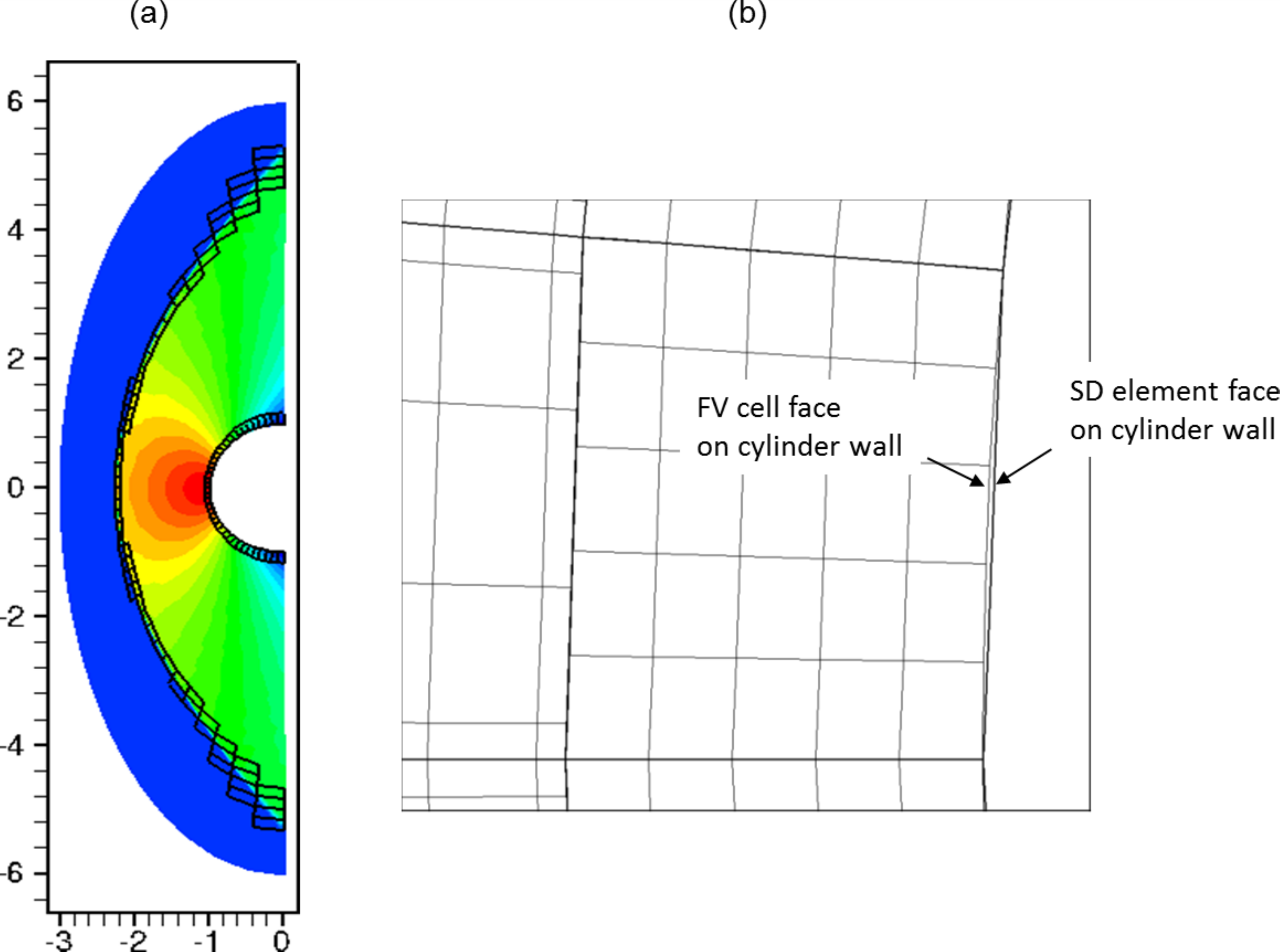} 
 \caption{Mach 2 flow over cylinder. (a) The locations of embedded FV elements indicated by black squares. (b) The close-up view of the embedded FV element on the cylinder wall.}
 \label{shock_cylinder_grid}
\end{center}
\end{figure}

\section{Discussion and Conclusion}
\label{conclusion}
In order to remedy the Gibbs phenomena in the spectral difference (SD) method associated with discontinuities in the flow fields, a hybrid SD/Embedded FV method for the hexahedral element is introduced.  In this proposed hybrid method, the FV element, consisting of structured FV subcells, is dynamically embedded in the base hexahedral elements which contain  discontinuities in the flow fields. The $5^{th}$-order WENO scheme with characteristic decomposition is currently employed as the shock-capturing scheme solving the governing equations on the embedded FV subcells. Away from the discontinuities, the $5^{th}$-order SD method is used in the smooth region.  The coupling between the SD element and the FV element is achieved  by the globally conserved mortar method.

In addition to the numerical error associated with the numerical schemes employed, the order of accuracy study shows that there are two additional sources of numerical error in the hybrid method, which are (a) the interpolation error in the interpolation of solutions between the set of Gauss points and  the set of FV subcell solution points in the SD element neighboring the FV element and (b) the interpolation error in the mortar projection process, interpolating the conservative variables from the FV subcells to the mortar interface and  interpolating the computed common fluxes back to the FV subcells.  It is also observed in the translating vortex test case that the numerical dissipation terms in mass and pressure fluxes in the AUSM$^+-$up scheme may introduce excessive numerical dissipation for unsteady low Mach number flows, especially when coupled with the high-order reconstruction scheme. The issue with the numerical dissipation associated with the AUSM$^+-$up scheme when coupled with the high-order reconstruction scheme needs further investigation, and it will be addressed in the future paper.  

Although the order of accuracy is degraded in the hybrid method, it is seen that the error norm in the hybrid method is still comparable to the error norm in the currently implemented $5^{th}$-order WENO scheme, e.g. in the case of translating density sine wave.   It should be noted that the purpose of embedding FV element with the WENO scheme is to capture the discontinuities, where the order of accuracy is essentially degraded in the WENO scheme, while  the solution in the smooth region is obtained by the less dissipative SD method. Thus, it is seen that the overall solution is  less affected by the reduction in the order of accuracy. This is demonstrated in 1D and 2D test cases comparing the simulations by the hybrid method and the standalone WENO scheme.
 
For the 1D and 2D test cases employed in this study,  overall the hybrid method gives comparable or better results compared with the standalone WENO scheme considering the same DOF of solution points in both methods. Particularly, when the complex flow structure dynamically develops, the hybrid method performs better compared with the standalone WENO scheme due to the low numerical dissipation in SD method. This is shown in the double Mach reflection test, where developing complex jet structure is better resolved with the hybrid method for the same solution DOF as the standalone WENO scheme.  In terms of grid size, hybrid method with much coarser grid is able to produce similar or better results as the ones on the finer grid with the standalone WENO scheme.  This may provide  benefits of applying the proposed hybrid method for complex geometry with easing the grid preparation efforts.  

The proposed hybrid method is presented in this paper using the structured hexahedral elements in order to deliver the main idea.  However, it should be mentioned that the proposed hybrid method is designed to be applied for the unstructured hexahedral elements.  Currently, this hybrid approach is being implemented for the unstructured hexahedral elements, and the results will be reported in the future paper together with the stability analysis including dispersion and dissipation properties, as discussed in \cite{martin2006, pirozzoli2002, kannan2011, kannan2012}.  In the case that the large number of FV subcells is used in the embedded FV element, such as the optimal number of FV subcells discussed in Dumbser {\it et al.} \cite{dumbser2014}, the discontinuity may be captured over a few subcells with the smooth or turbulent flows still covered by the majority of FV subcells. In such cases, low dissipative FV based shock-capturing scheme may be more appropriate to use in the embedded FV element.  Therefore, variants of WENO schemes with less numerical dissipation (e.g.,  MPWENO \cite{balsara2000}, WENO-Z \citep{borges2008}) or other high-order shock-capturing schemes (e.g., filter scheme \cite{yee2007}) with varying the order of the schemes and the number of FV subcells  will be examined and compared to study their effects on the numerical solutions and accuracy. The use of different Riemann solvers \cite{toro1999} with the high-order reconstruction scheme will also be investigated and reported in the future paper.

\noindent\\
\newblock{\bf Acknowledgement} 

The author would like to acknowledge Prof. Elaine S. Oran and Dr. Ryan W. Houim at  the University of Maryland for their valuable discussions and comments toward this work. In addition, the author acknowledges the University of Maryland supercomputing resources (http://www.it.umd.edu/hpcc) made available in conducting the research reported in this paper.       



\noindent\\
\newblock

\end{document}